\documentclass[10pt]{article}

\usepackage{amssymb,amsmath,amsfonts}
\usepackage{dsfont}
\usepackage{amsthm}
\usepackage{graphicx,color,subfig}
\usepackage{xspace}
\usepackage{mathrsfs}
\usepackage{pifont}

\newcommand{\Con}{\ensuremath{\mathscr C}}
\newcommand{\mb}[1]{\ensuremath{\mathbb{#1}}}
\newcommand{\N}{{\mb{N}}}
\DeclareMathOperator{\supp}{supp}
\newcommand{\R}{{\mb{R}}}
\DeclareMathSymbol{\intopmod}{\mathop}{symbols}{115}
\def\intmod{\intopmod}

\newtheorem{lemma}{Lemma}[section]
\newtheorem{theorem}[lemma]{Theorem}
\newtheorem{proposition}[lemma]{Proposition}
\newtheorem{corollary}[lemma]{Corollary}

\theoremstyle{definition}

\newtheorem{definition}[lemma]{Definition}
\newtheorem{remark}[lemma]{Remarks}
\newtheorem{case}{Case}
\newenvironment{prooff}[1][Proof]{\textbf{#1.} }{\hfill \rule{0.5em}{0.5em}}

\newcommand{\RR}{{\mathbb{R}}}

\newcommand{\dsp}{\displaystyle}

\makeatletter
\def\keywords{
    \vspace{1ex}
    \noindent
    \if@twocolumn
      \small{\bf  Keywords}\/---$\!$    \else
      \begin{center}\small\ {\bf Keywords}\end{center}\quotation\small
    \fi}
\def\endkeywords{\vspace{0.6em}\par\if@twocolumn\else\endquotation\fi
    \normalsize\rm}
\makeatother
\numberwithin{equation}{section}

\begin{document}

\title{The Kato Smoothing Effect for Regularized Schr\"odinger Equations in
Exterior Domains}
\author{Lassaad Aloui\thanks{Facult\'e des Sciences de Bizerte, email : lassaad.aloui@fsg.rnu.tn}, Moez
Khenissi\thanks{Universit\'e de Sousse, Ecole Sup. Sc. Tech. de Hammam Sousse, Tunisie,
email: moez.khenissi@fsg.rnu.tn}, Luc Robbiano\thanks{
Universit\'e de Versailles Saint-Quentin, Laboratoire de Math\'ematiques de
Versailles, CNRS UMR 8100, 45 Avenue des \'Etats-Unis, 78035 Versailles,
France, email : luc.robbiano@uvsq.fr}}
\maketitle

\begin{abstract}
We prove, under the exterior geometric control condition, the Kato smoothing
effect for solutions of an inhomogenous and damped Schr\"{o}dinger equation
on exterior domains.
\end{abstract}

\tableofcontents


\section{Introduction and results}

This paper is devoted to the study of a smoothing effect for a damped
Schr\"{o}dinger equation on exterior domain. In order to formulate the
results, we shall begin by recalling some results for Schr\"{o}dinger
equation linking the regularity of solutions and the geometry of domain
where these equations are posed. \newline
It is well known that the free Schr\"{o}dinger equation enjoys the property
of the $\Con^{\infty }$ smoothing effect, which can be described as follows:
For any distribution $u_{0}$ of compact support, the solution of the Cauchy
problem 
\begin{equation*}
\left\{ 
\begin{array}{l}
(i\partial _{t}+\Delta )u=0\text{ in }\mathbb{R}\times \mathbb{R}^{d} \\ [4pt] 
u_{|t=0}=u_{0},
\end{array}
\right.  
\end{equation*}
is infinitely differentiable with respect to $t$ and $x$ when $t\neq 0$ and 
$x\in \mathbb{R}^{d}$.

Another type of smoothing effect says that if 
$u_{0}\in L^{2}(\mathbb{R}^{d}) $ then the solution of the Schr\"{o}dinger equation satisfies the Kato 
$\frac{1}{2}$-smoothing effect ($H^{1/2}$-smoothing effect): 
\begin{equation*}
\int_{\mathbb{R}}\left\Vert \langle x\rangle ^{-s}\Delta ^{1/4}u\right\Vert
_{L^{2}(\mathbb{R}^{d})}^{2}\leq C\Vert u_{0}\Vert _{L^{2}}^{2},\text{\ }
s>1/2. 
\end{equation*}
This property of gain of regularity has been first observed in the case of 
$\mathbb{R}^{d}$ in the works of Constantin-Saut \cite{co.sa1}, Sj\"olin \cite{sjolin} and Vega~\cite{vega} and it has been extended 
locally in time
 to variable
coefficient operators with non trapping metric by Doi (\cite{doi1,Do})).

In the case of domains with boundary Burq, G\'erard and Tzvetkov \cite{b.g.t}
proved a local smoothing estimate for $\exp(it\Delta )$ in the exterior
domains with non-trapping assumption. Using the $TT^{\star }$ argument, the
proof of the smoothing effect with respect to initial data in \cite{b.g.t}
is reduced to the non-homogeneous bound which, by performing Fourier
transform in time, can be deduced from the bounds on the cut-off resolvent: 
\begin{equation*}
\Vert \chi (\lambda ^{2}-\Delta )^{-1}\chi \Vert _{L^{2}\rightarrow
L^{2}}\leq C,\text{ }\forall \lambda \gg 1.
\end{equation*}
The resolvent bound, for which the non-trapping assumption plays a crucial
role, is proven for $|\lambda |>>1$ in greater generality by Lax-Phillips ~\cite{LaxPh}, Melrose-Sjostrand~\cite{melJost,melJost2}, Vainberg~\cite{vainb} and Vazy-Zworski~\cite{V.Z} 

The Kato-effect has been extended by Robbiano and Zuily in \cite{RZ} to
variable coefficients operators with unbounded potential in exterior domains
with non trapping metric. The proof of their result is reduced to an
estimate localized in frequency which has been established by contradiction
using in a crucial way the semiclassical defect measure introduced by 
P.~Gerard~\cite{gerard} (see also \cite{Leb}). The use of the microlocal defect
measure to prove an estimate by contradiction method (Wilcox~\cite{wilcox})
go back\ to Lebeau~\cite{Leb}. This idea has been followed with success by
several authors (see Burq \cite{Bu,B2,bursmot} Aloui\ and
 Khenissi~\cite{alkh,alkh2,kh}).

In \cite{bursmot}, Burq proved that the non trapping condition is necessary
for the $H^{1/2}$ smoothing effect and showed, in the case of several convex
obstacles satisfying certain assumptions, the smoothing effect with an 
$\varepsilon >0$ loss: 
\begin{equation*}
\Vert \chi u\Vert _{L^{2}(H^{1/2-\varepsilon }(\Omega ))}\leq C\Vert
u_{0}\Vert _{L^{2}(\Omega )},
\end{equation*}
where $\chi$ is compactly supported.

On the other hand, the non-trapping assumption is also equivalent to the
uniform decay of the local energy for the wave equation 
(see \cite{LaxPh,rals,Melros}). For the trapping domains, when no such decay is hoped,
the idea of stabilization for the wave equation is to add a dissipative term
to the equation to force the energy of the solution to decrease uniformly.
There is a large literature on the problem of stabilization of wave
equation. In the case of bounded domains, we quote essentially the work of
J. Rauch and M. Taylor \cite{rauch} and the one of C. Bardos, G. Lebeau and
J. Rauch \cite{BLR} whose introduced and developed the geometric control
condition (GCC). This condition that asserts, roughly speaking, that every
ray of geometric optics enters the region where the damping term is
effective in a uniform time, turns out to be almost necessary and sufficient
for the uniform exponential decay of waves. In \cite{alkh}, Aloui and
Khenissi introduced the Exterior Geometric control condition (see below
Definition~\ref{egc})
 and hence extended the result of \cite{BLR} to the case of exterior
domains (see also \cite{alkh2} ).

Recently, by analogy with the stabilization problem the first author~\cite{al1,al2} has introduced the forced smoothing effect for Schr\"{o}dinger
equation in bounded domains; it consists to act on the equation to produce
some smoothing effects. More precisely he considered the following equation 
\begin{equation}
\left\{ 
\begin{array}{lll}
i\partial _{t}u-\Delta_{D} u +ia(x)(-\Delta_{D})^{\frac{1}{2}}a(x)u=0 & 
\text{in} & ]0,+\infty)\times \Omega, \\ 
u(0,.)=f & \text{in} & \Omega, \\ 
u|_{\mathbb{R}^{+}\times \partial \Omega }=0, &  & 
\end{array}
\right.  \label{eqr}
\end{equation}
where $\Omega$ is a bounded domain and $\Delta_{D}$ is the Dirichlet-Laplace
operator on $\Omega$.

Using the strategy of \cite{b.g.t}, Aloui~\cite{al2} proved a weak Kato
-Smoothing effect:

\begin{equation}
\left\Vert v\right\Vert _{L^{2}([\varepsilon ,T],H_{D}^{s+1}(\Omega ))}\leq
c\left\Vert v_{0}\right\Vert _{H_{D}^{s}(\Omega )},  \label{fai}
\end{equation}
where $0<\varepsilon <T<\infty $ and $v_{0}\in H_{D}^{s}(\Omega )$, (See 
\cite{al2} for the definition of $H_{D}^{s}$).

By iteration of the last result, Aloui deduced also a 
$\Con^{\infty }$-smoothing effect for the regularized Schr\"{o}dinger equation 
(\ref{eqr}).
Recently, Aloui, Khenissi and Vodev~\cite{alKhVo} have proved that the
Geometric control condition is not necessary to obtain the forced 
$\Con^{\infty }$- smoothing effect.

On the other hand, using the arguments of \cite{b.g.t}, we can prove, for
the equation (\ref{eqr}) in exterior domains, the cut-off resolvent bound,
which is sufficient to deduce the non-homogenous bound. But, unfortunately,
the generator operator $\Delta _{D}-ia(x)(-\Delta _{D})^{\frac{1}{2}}a(x)$
is not self-adjoint and then the $TT^{\star }$ argument fails. For this
reason, we can not prove (with this strategy) the weak Kato-smoothing effect
(\ref{fai}) for exterior domains.

The question now is the following:

Can we establish the Kato-smoothing effect for the regularized 
Schr\"{o}dinger equation (\ref{eqr}) for which the Geometric Control Condition is
necessary? and if so, does this result still hold for exterior problems?

In this paper, we give an affirmative answer. Indeed, under the Exterior
Geometric Control condition, we prove the Kato-smoothing effect and the non
homogenous bound for the regularized Schr\"{o}dinger equation in exterior
domains. Notice that the case of bounded domains can be treated by the same
method.

Our approach for deriving such results is to combine the strategies of
Robbiano-Zuily in \cite{RZ} \ and Aloui-Khenissi in \cite{alkh}, \cite{kh}.

In order to state our results, we give several notations and assumptions. 
\newline
Let $K$ be a compact obstacle in $\mathbb{R}^{d}$ whose complement $\Omega $
an open set with $\Con^{\infty }$ boundary $\partial \Omega $ and $\tilde{P}
$ be a second-order differential operator of the form 

\begin{equation}  \label{eq:P}
\tilde{P}=\sum_{j,k=1}^{d}D_j(b^{jk}D_k)+V(x), \qquad D_j=\frac{\partial}{
i\partial x_j},
\end{equation}

where coefficients $b^{jk}$\ and $V$\ are assumed to be in 
$\Con^{\infty }(\mathbb{R}^d),$\ real valued,\ and 
$b^{jk}=b^{kj},$ $1\leq j,$ $k\leq d.$

Throughout this paper, 
$\left\langle x\right\rangle :=(1+|x|^{2})^{\frac{1}{2}}$ 
and we denote by $S_{\Omega }(M,g)$ the H\"{o}rmander's class of symbols
if $M$ is a weight and the metric 
\begin{equation*}
g=\frac{dx^{2}}{\left\langle x\right\rangle ^{2}}+
\frac{d\xi ^{2}}{\left\langle \xi \right\rangle ^{2}}.
\end{equation*}
We shall denote by $p$ the principal symbol of $\tilde{P}$, namely 
\begin{equation*}
p(x,\xi )=\sum_{j,k=1}^{d}b^{jk}(x)\xi _{j}\xi _{k},
\end{equation*}
and we assume that 
\begin{equation}
\exists \text{ }c>0:p(x,\xi )\geq c|\xi |^{2},\text{ \ for }x\text{ in }
\mathbb{R}^{d}\text{\ and }\xi \text{\ in }\mathbb{R}^{d},  \label{eq:illip}
\end{equation}
\begin{equation}
\left\{ 
\begin{array}{l}
(i)\text{ }b^{jk}\in S_{\Omega }(1,g),\text{ }\nabla _{x}b^{jk}(x)=
o(\frac{1}{|x|}),\text{ \ }|x|\rightarrow +\infty ,\text{ \ }1\leq j,\text{ }k\leq d.
\\ 
(ii)\text{ }V\in S_{\Omega }(\left\langle x\right\rangle ^{2},g),\text{ \ }
V\geq -C_{0}\text{ for some positive constant }C_{0}.
\end{array}
\right.  \label{hyp1}
\end{equation}

Under the assumptions (\ref{eq:illip}) and (\ref{hyp1}), the operator 
$\tilde{P}$ is essentially self-adjoint on $\Con^\infty_0(\Omega)$ and we
denote by $P$ its self-adjoint extension. \newline
Now we set 
\begin{equation*}
\Lambda =((1+C_{0})Id+P)^{1/2},
\end{equation*}
which is well defined by functional calculus of self-adjoint positive
operators. \newline
We consider the following regularized Schr\"{o}dinger equation
\begin{equation}
\left\{ 
\begin{array}{l}
(D_{t}+P)u-ia\Lambda au=f\text{ in }]0,+\infty)\times \Omega \\[4pt] 
u=0\text{ on }[0,+\infty)\times \partial \Omega , \\ 
u_{|t=0}=u_{0},
\end{array}
\right.  \label{eq: Equa}
\end{equation}
where $(u_{0},f)\in \Con_{0}^{\infty }(\Omega)\times\Con_{0}^{\infty }(]0,+\infty)\times\Omega )$ and 
$a\in \Con_{0}^{\infty }(\overline{\Omega }).$ 

Let's recall the Exterior Geometric Control (E.G.C.) condition \cite{alkh}

\begin{definition}[E.G.C.]
\label{egc}Let $R>0$ be such that $K\subset B_{R}=\{|x|<R\}$ and $\omega $ be a subset
of $\Omega .$ We say that $\omega $ verifies the Exterior Geometric Control
condition on $B_{R}$ (E.G.C.) if there exists $T_{R}>0$\ such that every
generalized bicharacteristic $\gamma $ starting from $B_{R}$ at time $t=0,$
is such that:

\begin{itemize}
\item $\gamma $ leaves $\mathbb{R}^{+}\times B_{R}$ before the time $T_{R},$
or

\item $\gamma $ meets $\mathbb{R}^{+}\times \omega $ between the times $0$
and $T_{R}.$
\end{itemize}
\end{definition}

We assume also that the bicaracteristics have no contact of infinite order
with the boundary (see, for a precise statement, Definition~\ref{contact
d'ordre infini}).

Under this condition on $\omega =\{x\in \Omega ,a^{2}(x)>0\},$ we can state
our main result.

\begin{theorem} 
\label{A}Let $T>0$, $\alpha \in (-1/2,1/2)$ and $s\in (1/2,1]$. Let $P$
defined by (\ref{eq:P}) satisfying the assumptions (\ref{eq:illip}) and (\ref
{hyp1}). Then under, the E.G.C on $\omega$ one can find a positive constant $C(T,\alpha
,s)=C$ such that 
\begin{equation}
\int_{0}^{T}\left\Vert \Lambda ^{\alpha +1/2}\langle x\rangle
^{-s}u\right\Vert _{L^{2}(\Omega )}^{2}dt+\!\sup_{t\in \lbrack 0,T]}\Vert
\Lambda ^{\alpha }u(t)\Vert _{L^{2}(\Omega )}^{2}\!\leq \!C\left( \Vert
\Lambda ^{\alpha }u_{0}\Vert _{L^{2}(\Omega )}^{2}+\int_{0}^{T}\left\Vert
\Lambda ^{\alpha -1/2}\langle x\rangle ^{s}f\right\Vert _{L^{2}(\Omega
)}^{2}dt\right)  \label{eq:estmGlob}
\end{equation}
for all $u_{0}$ 
in $\Con_{0}^{\infty }(\Omega )$, 
$f$ in $\Con_{0}^{\infty }(\Omega \times \mathbb{R}^{+} )$, where $u$ denotes the
solution of (\ref{eq: Equa}).
\end{theorem}

Working with $\tilde{u}=e^{i(1+C_{0})t}u,$ one may assume $V\geq 1$ in 
(\ref{hyp1}) and $\Lambda =P^{1/2},$ which will be assumed in the sequel. It
turns into the following equation 
\begin{equation}
\left\{ 
\begin{array}{l}
(D_{t}+P)u-iaP^{1/2}au=f\text{ in }[0, +\infty)\times \Omega \\[4pt] 
u=0\text{ on }[0,\infty)\times \partial \Omega , \\ 
u_{|t=0}=u_{0},
\end{array}
\right.
\end{equation}
where $P\geq 1.$

\begin{remark}
$
\begin{array}{c}
\end{array}
$

\begin{enumerate}
\item When the obstacle is nontrapping, we obtain the result of Robbiano
Zuily \cite{RZ} by\ taking $a(x)=0$ and moreover, we improve their result to
non homogenous bound.

\item If we consider the equation in a bounded domain $\Omega $ of $\mathbb{R
}^{d},$ and replace the exterior geometric condition (E.G.C) by the
classical microlocal condition of Bardos-Lebeau-Rauch \cite{BLR}, we can
still prove the Kato-effect and then we improve the result of 
Aloui~\cite{al2}.

\item If there is a trapped ray which does not intersect the regularized
region, due to Burq \cite{bursmot}, the Kato-effect does not hold. In this
context, our result is thus optimal.
\end{enumerate}
\end{remark}

The rest of the paper is organized as follows: Section~\ref{proofs} is devoted to the
proof of Theorem \ref{A} while in the Section~\ref{appendix} we shall prove some Lemmata
used in Section~\ref{proofs}.

\section{Proofs}\label{proofs}

 Let's describe the strategy of  the proof of theorem 1.2. In a first step, we reduce the estimate \eqref{eq:estmGlob} to an analogue one localized in frequencies. By following a contradiction argument, we can construct an adapted microlocal defect measure. Our aim in the rest of the proof is to obtain a contradiction on this measure. First, we prove that this measure is not identically null.  Next, we show that it is null on incoming set and on $\{ a^2>0\}$. Finaly, using the geometrical assumption (E.G.C.) and that the support of this measure is propagated along the generalized flow, we conclude that the measure is identically null. This gives the contradiction.

\subsection{Reduction to an estimate localized in frequency}

\label{reduction en frequences}

We recall the Paley-Littlewood decomposition. Let $\Phi \in \Con_{0}^{\infty
}([0,+\infty ))$ be a decreasing function such that
\begin{equation*}
\Phi (s)=1\text{ if \ }s\leq 1/2,\text{ \ }\Phi (s)=0\text{ \ if \ }s\geq 1.
\newline
\end{equation*}
Let $\psi (s)=\Phi (4^{-1}s)-\Phi (s)$, $\psi (s)=0$ if $s\leq 1/2$ or 
$s\geq 4$, $0\leq \psi \leq 1$. For $s\geq 0$ we have 
\begin{equation*}
\displaystyle1=\Phi (s)+\sum_{n=0}^{+\infty }\psi (4^{-n}s),
\end{equation*}
and using $P\geq 1$, we have 
\begin{equation*}
\displaystyle u=\sum_{n=0}^{+\infty }\psi (4^{-n}P)u.
\end{equation*}
For support reason 
\begin{equation*}
\psi (4^{-n}s)\psi (4^{-k}s)=0\text{ if }|k-n|\geq 2,
\end{equation*}
thus there exists $C>0$ such that for all $u\in L^{2}(\Omega )$, 
\begin{equation*}
\Vert u\Vert _{L^{2}(\Omega )}^{2}\leq C\sum_{n=0}^{+\infty }\Vert \psi
(4^{-n}P)u\Vert _{L^{2}(\Omega )}^{2}\leq C^{2}\Vert u\Vert _{L^{2}(\Omega
)}^{2}. 
\end{equation*}
In the sequel we denote by $h_{n}=2^{-n}$ and $u_{n}=u_{h_{n}}=\psi
(h_{n}^{2}P)u $.\newline
If $u$ satisfies 
\begin{equation}
D_{t}u+Pu-iaP^{1/2}(au)=f,  \label{eq:Schrod}
\end{equation}
thus $u_{n}$ is a solution of the following semi-classical Schr\"odinger equation:  
\begin{equation}
h_{n}^{2}(D_{t}+P)u_{n}-ih_{n}a(h_{n}^{2}P)^{1/2}(au_{n})=h_{n}g_{n},
\label{eq:lowfreq}
\end{equation}
where 
\begin{equation}
g_{n}=g_{h_{n}}=h_{n}\psi (h_{n}^{2}P)f+i[\psi
(h_{n}^{2}P),a](h_{n}^{2}P)^{1/2}(au)+ia(h_{n}^{2}P)^{1/2}[\psi
(h_{n}^{2}P),a]u.  \label{eq:g}
\end{equation}

\begin{proposition}
\label{prop:lowfreq} Let $s\in (1/2,1]$, $T>0$ and $\alpha \in (-1/2,1/2)$.
Assume there exists $C>0$ such that for $u_{n}=\psi (h_{n}^{2}P)u$
satisfying \eqref{eq:lowfreq}, we have, for all $n\geq 1$
\begin{equation}
\Vert \langle x\rangle ^{-s}u_{n}\Vert _{L^{2}([0,T]\times \Omega
)}^{2}+h_{n}\sup_{t\in \lbrack 0,T]}\Vert u_{n}(t)\Vert _{L^{2}(\Omega
)}^{2}\leq C\left( h_{n}\Vert u_{n}(0)\Vert _{L^{2}(\Omega )}^{2}+\Vert
\langle x\rangle ^{s}g_{n}\Vert _{L^{2}([0,T]\times \Omega )}^{2}\right),
\label{prop:lowfreq1}
\end{equation}
then there exists $C^{\prime }>0$ such that for all $u$ satisfying 
\eqref{eq:Schrod} we have 
\begin{equation}  \label{inegalite avec P alpha}
\begin{split}
\| P^{\alpha/2+1/4}\langle x\rangle^{-s}u\|_{L^2([0,T]\times
\Omega)}^2+\sup_{t\in[0,T]}\| P^{\alpha/2} u(t)\|_{L^2(\Omega)}^2 \qquad
\qquad\qquad \qquad \qquad \qquad\qquad \\
\le C^{\prime } \left( \| P^{\alpha/2}u(0)\|_{L^2(\Omega)}^2+\|
P^{\alpha/2-1/4}\langle x\rangle^s f\|_{L^2([0,T]\times \Omega)}^2 \right) .
\end{split}
\end{equation}
\end{proposition}

\begin{prooff}
We multiply \eqref{prop:lowfreq1} by $h_{n}^{-2\alpha -1}$ and we sum over 
$n\in \N$, we obtain, 
\begin{equation}
\begin{split}
& \sum_{n\in \N}h_{n}^{-2\alpha -1}\Vert \langle x\rangle ^{-s}u_{n}\Vert
_{L^{2}([0,T]\times \Omega )}^{2}+\sum_{n\in \N}h_{n}^{-2\alpha }\sup_{t\in
\lbrack 0,T]}\Vert u_{n}(t)\Vert _{L^{2}(\Omega )}^{2}\qquad \qquad \\
& \leq C\left( \sum_{n\in \N}h_{n}^{-2\alpha }\Vert u_{n}(0)\Vert
_{L^{2}(\Omega )}^{2}+\sum_{n\in \N}h_{n}^{-2\alpha -1}\Vert \langle
x\rangle ^{s}g_{n}\Vert _{L^{2}([0,T]\times \Omega )}^{2}\right) .
\end{split}
\label{inegalite sur la somme en n}
\end{equation}
Now, let us estimate each term appearing in inequality 
\eqref{inegalite
avec P alpha}. We have, 
\begin{align}
\sup_{t\in \lbrack 0,T]}\Vert P^{\alpha /2}u(t)\Vert _{L^{2}(\Omega )}^{2}&
\leq C\sup_{t\in \lbrack 0,T]}\sum_{n\in \N}\Vert \psi (h_{n}^{2}P)P^{\alpha
/2}u(t)\Vert _{L^{2}(\Omega )}^{2}  \notag \\
& \leq C\sup_{t\in \lbrack 0,T]}\sum_{n\in \N}h_{n}^{-2\alpha }\Vert \psi
_{0}(h_{n}^{2}P)u(t)\Vert _{L^{2}(\Omega )}^{2}\text{ where }\psi
_{0}(\sigma )=\sigma ^{\alpha /2}\psi (\sigma )  \notag \\
& \leq C\sum_{n\in \N}h_{n}^{-2\alpha }\sup_{t\in \lbrack 0,T]}\Vert \psi
(h_{n}^{2}P)u(t)\Vert _{L^{2}(\Omega )}^{2}.
\label{Premiere inegalite P alpha}
\end{align}
We have also with 
$\psi _{1}(\sigma )=\sigma ^{\alpha /2+1/4}\psi (\sigma )$
, 
\begin{align}
\Vert P^{\alpha /2+1/4}\langle x\rangle ^{-s}u\Vert _{L^{2}([0,T]\times
\Omega )}^{2}& \leq C\sum_{n\in \N}h_{n}^{-2\alpha -1}\Vert \psi
_{1}(h_{n}^{2}P)\langle x\rangle ^{-s}u\Vert _{L^{2}([0,T]\times \Omega
)}^{2}  \notag \\
& \leq C\sum_{n\in \N}h_{n}^{-2\alpha -1}\Vert \langle x\rangle ^{-s}\psi
(h_{n}^{2}P)u\Vert _{L^{2}([0,T]\times \Omega )}^{2}
\text{( by Lemma~\ref{equivalence norme H alpha} )}  \notag \\
& \leq C\sum_{n\in \N}h_{n}^{-2\alpha -1}\Vert \langle x\rangle
^{-s}u_{n}\Vert _{L^{2}([0,T]\times \Omega )}^{2}.
\label{deuxieme
inegalite P alpha}
\end{align}
Now we can estimate, with $\psi _{2}(\sigma )=\sigma ^{-\alpha /2}\psi
(\sigma )$, 
\begin{align}
\sum_{n\in \N}h_{n}^{-2\alpha }\Vert u_{n}(0)\Vert _{L^{2}(\Omega )}^{2}&
\leq C\sum_{n\in \N}\Vert \psi _{2}(h_{n}^{2}P)P^{\alpha /2}u(0)\Vert
_{L^{2}(\Omega )}^{2}  \notag \\
& \leq C\Vert P^{\alpha /2}u(0)\Vert _{L^{2}(\Omega )}^{2}.
\label{troisieme
inegalite P alpha}
\end{align}
The term $g_{n}$ contains three terms (see \eqref{eq:g}). For the first, we
have, with $\psi _{3}(\sigma )=\sigma ^{-\alpha /2+1/4}\psi (\sigma )$, 
\begin{align}
\sum_{n\in \N}h_{n}^{-2\alpha +1}\Vert \langle x\rangle ^{s}\psi
(h_{n}^{2}P)f\Vert ^2_{L^{2}([0,T]
\times \Omega )}& \leq \sum_{n\in \N
}h_{n}^{-2\alpha +1}\Vert \psi (h_{n}^{2}P)\langle x\rangle ^{s}f\Vert^2
_{L^{2}([0,T]
\times \Omega )}  \notag \\
& \leq C\sum_{n\in \N}\Vert \psi _{3}(h_{n}^{2}P)P^{\alpha /2-1/4}\langle
x\rangle ^{s}f\Vert _{L^{2}([0,T]\times \Omega )}^{2}  \notag \\
& \leq C\Vert P^{\alpha /2-1/4}\langle x\rangle ^{s}f\Vert
_{L^{2}([0,T]\times \Omega )}^{2}.  \label{estimation premier term
g_n}
\end{align}
For the second and the third terms of $g_{n}$ we can apply the Lemmata~\ref
{lemma : premier commutateur} and \ref{lemme deuxieme terme}, to obtain with 
\eqref{estimation premier term g_n}, 
\begin{equation}
\sum_{n\in \N}h_{n}^{-2\alpha -1}\Vert \langle x\rangle ^{s}g_{n}\Vert
_{L^{2}([0,T]\times \Omega )}^{2}\leq C\Vert P^{\alpha /2-1/4}\langle
x\rangle ^{s}f\Vert _{L^{2}([0,T]\times \Omega )}^{2}+C\Vert P^{\alpha
/2}u\Vert _{L^{2}([0,T]\times \Omega )}^{2}.
\label{Quatrieme inegalite P alpha}
\end{equation}
Then following \eqref{inegalite sur la somme en n} 
\eqref{Premiere
inegalite P alpha}, \eqref{deuxieme inegalite P alpha}, 
\eqref{troisieme
inegalite P alpha} and \eqref{Quatrieme inegalite
P alpha}, we obtain 
\begin{equation*}
\begin{split}
& \Vert P^{\alpha /2+1/4}\langle x\rangle ^{-s}u\Vert _{L^{2}([0,T]\times
\Omega )}^{2}+\sup_{t\in \lbrack 0,T]}\Vert P^{\alpha /2}u(t)\Vert
_{L^{2}(\Omega )}^{2}\qquad \qquad \qquad \qquad \qquad \qquad \qquad \\
& \leq C\left( \Vert P^{\alpha /2}u(0)\Vert _{L^{2}(\Omega )}^{2}+\Vert
P^{\alpha /2-1/4}\langle x\rangle ^{s}f\Vert _{L^{2}([0,T]\times \Omega
)}^{2}+\Vert P^{\alpha /2}u\Vert _{L^{2}([0,T]\times \Omega )}^{2}\right) .
\end{split}
\end{equation*}
By Gronwall's Lemma, we can remove the last term in the previous inequality
and we obtain~\eqref{inegalite avec P alpha}.

\end{prooff}

\subsection{Construction of microlocal defect measure}

In this section we will prove the localized frequency 
estimate~(\ref{prop:lowfreq1}) by a contradiction argument and using microlocal defect
measure.

More precisely, let $u_{h}$ solution of 
\begin{equation}
h^{2}(D_{t}+P)u_{h}-iha(h^{2}P)^{1/2}(au_{h})=hg_{h}.  \label{eq:lowfre2}
\end{equation}
We will prove by contradiction the following estimate, 
\begin{equation}
\Vert \langle x\rangle ^{-s}u_{h}\Vert _{L^{2}([0,T]\times \Omega
)}^{2}+h\sup_{t\in \lbrack 0,T]}\Vert u_{h}(t)\Vert _{L^{2}(\Omega
)}^{2}\leq Ch\Vert u_{h}(0)\Vert _{L^{2}(\Omega )}^{2}+C\Vert \langle
x\rangle ^{s}g_{h}\Vert _{L^{2}([0,T]\times \Omega )}^{2}.
\label{eq:lowfreqIstim2}
\end{equation}

Assuming it is false. Taking $C=k\in \mathbb{N}$, we deduce sequences 
$h_{k}\underset{k\rightarrow +\infty }{\rightarrow }0,$ 
$u_{k}^{0}=u_{h_k}(0)\in
L^{2}(\Omega )$ and $g_{k}=g_{h_k}\in L^{2}(\Omega )$ such that,
\begin{equation}
h_{k}\left\Vert u_{k}^{0}\right\Vert _{L^{2}(\Omega )}^{2}
\underset{k\rightarrow +\infty }{\rightarrow }0,\text{ }\left\Vert \left\langle
x\right\rangle ^{s}g_{k}\right\Vert _{L^{2}([0,T]\times \Omega )}^{2}
\underset{k\rightarrow +\infty }{\rightarrow }0.  \label{eq:contradiction}
\end{equation}
We normalize by the left term in (\ref{eq:lowfreqIstim2}), thus 
\begin{equation*}
\left\Vert \left\langle x\right\rangle ^{-s}u_{k}\right\Vert
_{L^{2}([0,T]\times \Omega )}^{2}+\text{ }h_{k}\sup_{t\in \lbrack
0,T]}\left\Vert u_{k}(t)\right\Vert _{L^{2}(\Omega )}^{2}=1,
\end{equation*}
where, for simplicity, we have denoted $u_{h_k}=u_k$. By the 
Lemma~\ref{lemma:A} we have 
\begin{equation}
\text{ }h_{k}\sup_{t\in \lbrack 0,T]}\left\Vert u_{k}(t)\right\Vert
_{L^{2}(\Omega )}^{2}\underset{k\rightarrow +\infty }{\rightarrow }0,
\label{eq:contradiction2}
\end{equation}
then 
\begin{equation}
\left\Vert \left\langle x\right\rangle ^{-s}u_{k}\right\Vert
_{L^{2}([0,T]\times \Omega )}^{2}\underset{k\rightarrow +\infty }
{\rightarrow }1 .  \label{eq:contradiction3}
\end{equation}
The sequence $(u_{k})$ is bounded in $L_{loc}^{2}
(\mathbb{R}_{t},L_{loc}^{2}(\Omega )).$ Indeed, for $R>0$ , there exists $c>0$ such
that $\left\langle x\right\rangle ^{-2s}\geq c,\,\;\forall x\in B(0,R)$ and then we have
\begin{equation}
\int_{0}^{T}\int_{\Omega \cap B_{R}}|u_{k}|^{2}dtdx\leq \frac{1}{c}
\int_{0}^{T}\int_{\Omega \cap B_{R}}\left\langle x\right\rangle
^{-2s}|u_{k}|^{2}dtdx\leq \frac{1}{c}.  \label{eq:bound}
\end{equation}
We set 
\begin{equation}
\left\{ 
\begin{array}{c}
w_{k}=1_{\Omega }u_{k}(t) \\ 
W_{k}=1_{[0,T]}w_{k}.
\end{array}
\right.  \label{eq:plong}
\end{equation}
It follows from (\ref{eq:bound}) that the sequence $(W_{k})$ is bounded in 
$L^{2}(\mathbb{R}_{t},L_{loc}^{2}(\mathbb{R}^{d})).$\newline
We associate to a symbol $b=b(x,t,\xi ,\tau )\in \Con_{0}^{\infty }(T^{\ast }
{\mathbb{R}}^{d+1})$ the semiclassical pseudo-differential operator (pdo) by
the formula 
\begin{equation*}
{\mathcal{O}}p(b)(y,s,hD_{x},h^{2}D_{t})v(x,t)=\frac{1}{(2\pi h)^{d+1}}\iint
e^{i\left( \frac{x-y}{h}\xi +\frac{t-s}{h^{2}}\tau \right) }\varphi
(y)b(x,t,\xi ,\tau )v(y,s)dydsd\xi d\tau,
\end{equation*}
where $\varphi \in \Con_{0}^{\infty }({\mathbb{R}}^{d})$ is equal to one on
a neighborhood of the $x$-projection of the support of $b$. As in \cite{RZ}
we can associate to $(W_{k})$ a semi-classical measure $\mu .$ More
precisely,

\begin{proposition}
\label{mesure}There exists a subsequence $(W_{\sigma (k)})$ and a Radon
measure $\mu $ on $T^{\ast }\mathbb{R}^{d+1}$ such that for every 
$b\in \Con_{0}^{\infty }(T^{\ast }{\mathbb{R}}^{d+1})$ one has 
\begin{equation*}
\lim_{k\rightarrow +\infty }\left( \mathcal{O}p(b)\left( x,t,h_{\sigma
(k)}D_{x},h_{\sigma (k)}^{2}D_{t}\right) W_{\sigma (k)},W_{\sigma
(k)}\right) _{L^{2}({\mathbb{R}}^{d+1})}=\left\langle \mu ,b\right\rangle .
\end{equation*}
\end{proposition}

We prove first that the measure $\mu $ satisfies the  
  following property.
\begin{proposition}
\label{prop:sopport}The support of $\mu $ is contained in the characteristic
set of the operator $D_{t}+P$ 
\begin{equation}
\Sigma =\{(x,t,\xi ,\tau )\in T^{\ast }\mathbb{R}^{d+1}:x\in
 \overline{\Omega },t\in \lbrack 0,T]\text{ and }\tau +p(x,\xi )=0\}.
\label{prop:carac}
\end{equation}
\end{proposition}

\begin{prooff}
According to (\ref{eq:plong}), it is obvious that 
\begin{equation*}
\supp\mu \subset \{(x,t,\xi ,\tau )\in T^{\ast }\mathbb{R}^{d+1}:x\in 
\overline{\Omega },t\in \lbrack 0,T]\}.
\end{equation*}
Therefore it remains to show that if $m_{0}=(x_{0},t_{0},\xi _{0},\tau _{0})$
with $x_{0}\in \overline{\Omega },t_{0}\in \lbrack 0,T],$ and $\tau
_{0}+p(x_{0},\xi _{0})\neq 0$ then $m_{0}\notin \supp \mu .$ For simplicity, we 
shall denote
the sequence $W_{\sigma(k)}$ by $W_k$.
\begin{case}
Assume that $x_{0}\in \Omega .$

Let $\varepsilon >0$ be such that $B(x_{0},\varepsilon )\subset \Omega $,
 $\varphi \in \Con_{0}^{\infty }(B(x_{0},\varepsilon )),$ $\varphi =1$ on 
$B(x_{0},\frac{\varepsilon }{2})$ and $\tilde{\varphi}\in \Con_{0}^{\infty
}(\Omega ),$ $\tilde{\varphi}=1$ on $\supp\varphi .$ Let $b\in 
\Con_{0}^{\infty }(\mathbb{R}_{x}^{d}\times \mathbb{R}_{\xi }^{d})$ such that 
$\pi _{x}\supp b\subset B(x_{0},\frac{\varepsilon }{2})$ and 
$\chi \in \Con_{0}^{\infty }(\mathbb{R}_{t}\times \mathbb{R}_{\tau }).$ Recall that we
have $W_{k}=1_{[0,T]}1_{\Omega }u_{k}$ and that $(u_{k})$ is bounded
sequence in $L^{2}([0,T],L_{loc}^{2}(\Omega )).$ We set 
\begin{equation*}
I_{k}=\left( b(x,h_{k}D_{x}\right) \chi (t,h_{k}^{2}D_{t})\varphi
(x)h_{k}^{2}(D_{t}+P(x,D_{x}))W_{k},\tilde{\varphi}W_{k})_{L^{2}
(\mathbb{R}^{d+1})}.
\end{equation*}
As in \cite{RZ} we have 
\begin{equation}
\lim_{k\rightarrow +\infty }I_{k}=\left\langle \mu ,(\tau +p)b\chi
\right\rangle .  \label{eq:carac}
\end{equation}
On the other hand, since we have 
\begin{equation*}
h_{k}^{2}(D_{t}+P(x,D_{x}))u_{k}=h_{k}ia(h_{k}^{2}P)^{1/2}au_{k}+h_{k}g_{k},
\end{equation*}
and $\varphi \in \Con_{0}^{\infty }(\Omega )$,
\begin{equation}
\varphi (h_{k}^{2}D_{t}+h_{k}^{2}P(x,D_{x}))W_{k}=\varphi
(ih_{k}a(h_{k}^{2}P)^{1/2}au_{k}+h_{k}g_{k})+h_k^2\varphi (u_{k}(0)\delta
_{t=0}-h_k^2u_{k}(T)\delta _{t=T}) .  \label{eq:terms}
\end{equation}
Then $I_{k}$ is a sum of four terms,
\begin{equation*}
I_{k}=I_{k}^{1}+I_{k}^{2}+I_{k}^{3}+I_{k}^{4},
\end{equation*}

\begin{align*}
&I_{k}^{1} =ih_{k}\left( b(x,h_{k}D_{x}\right) \chi
(t,h_{k}^{2}D_{t})\varphi (x)a(h_{k}^{2}P)^{1/2}au_{k},\tilde{\varphi}
W_{k})_{L^{2}(\mathbb{R}^{d+1})} \\
&I_{k}^{2} =h_{k}\left( b(x,h_{k}D_{x}\right) \chi (t,h_{k}^{2}D_{t})\varphi
(x)g_{k},\tilde{\varphi}W_{k})_{L^{2}(\mathbb{R}^{d+1})} \\
&I_{k}^{3} =\left( b(x,h_{k}D_{x}\right) \chi (t,h_{k}^{2}D_{t})h_k^2\varphi
(x)u_{k}(0)\delta _{t=0},\tilde{\varphi}W_{k})_{L^{2}(\mathbb{R}^{d+1})} \\
&I_{k}^{4} =-(b(x,h_{k}D_{x})\chi (t,h_{k}^{2}D_{t})h_k^2\varphi
(x)u_{k}(T)\delta _{t=T},\tilde{\varphi}W_{k})_{L^{2}(\mathbb{R}^{d+1})}.
\end{align*}
For the first term $I_{k}^{1}$, we use the Lemma~\ref{lemma:H}, 
we have, 
\begin{equation}
\left\Vert (h_{k}^{2}P)^{1/2}au_{k}\right\Vert _{L^{2}(\Omega )}^2\leq
Ch_{k}^{2}\Vert u_{k}\Vert _{L^{2}(\Omega )}^{2}+C\Vert au_{k}\Vert
_{L^{2}(\Omega )}^{2} ,  \label{eq:auk}
\end{equation}
and we deduce, 
\begin{equation}
|I_{k}^{1}|\leq c(h_{k}^{2}\sup_{t\in \lbrack 0,T]}\Vert u_{k}\Vert
_{L^{2}(\Omega )}^{2}+h_{k}\sup_{t\in \lbrack 0,T]}\Vert u_{k}\Vert
_{L^{2}(\Omega )}^{2}).  \label{eq:I11}
\end{equation}
Then we obtain, that $I^1_k$ goes to zero by \eqref{eq:contradiction2}. 
For the second term $I_{k}^{2}$, 
\begin{align*}
\left\vert I_{k}^{2}\right\vert & \leq h_{k}\left\Vert g_{k}\right\Vert
_{L^{2}([0,T],B(x_{0},\varepsilon ))}\left\Vert \tilde{\varphi}
W_{k}\right\Vert _{L^{2}(\mathbb{R}^{d+1})} \\
& \leq Ch_{k}\left\Vert \left\langle x\right\rangle ^{s}g_{k}\right\Vert
_{L^{2}([0,T]\times \Omega )}\left\Vert \left\langle x\right\rangle
^{-s}u_{k}\right\Vert _{L^{2}([0,T]\times \Omega )}.
\end{align*}
Using (\ref{eq:contradiction}) and (\ref{eq:contradiction3}), we deduce that
\begin{equation}
\lim_{k\rightarrow +\infty }I_{k}^{2}=0.  \label{eq:I2}
\end{equation}
The third and fourth terms in (\ref{eq:terms}) have the following form,
\begin{equation*}
J_{k}=\left( b(x,h_{k}D_{x})\chi (t,h_{k}^{2}D_{t})\varphi
h_{k}^{2}u_{k}(s)\delta _{t=s},\tilde{\varphi}W_{k}\right) _{L^{2}
(\mathbb{R}^{d+1})},\text{ \ \ }s=0\text{ or }T.
\end{equation*}
Since $(\tilde{\varphi}W_{k})$ is bounded in $L^{2}(\mathbb{R}^{d+1}),$ we
see that 
\begin{equation*}
|J_{k}|^{2}\leq c\left\Vert b\varphi w_{k}(s)\right\Vert _{L^{2}
(\mathbb{R}^{d})}^{2}\left\Vert h_{k}^{2}\chi (t,h_{k}^{2}D_{t})\delta
_{t=s}\right\Vert _{L^{2}(\mathbb{R})}^{2}\sup_{t\in \lbrack 0,T]}\Vert
u_{k}(t)\Vert _{L^{2}(\Omega )}^{2},
\end{equation*} 
so, using \cite[Lemma A.5]{RZ} with $p=2$ and $l=2$, we deduce that,
\begin{equation}
|J_{k}|^{2}\leq ch_{k}^{2}\left\Vert u_{k}(s)\right\Vert _{L^{2}(\Omega
)}^{2}\sup_{t\in \lbrack 0,T]}\Vert u_{k}(t)\Vert _{L^{2}(\Omega )}^{2}\leq
c\,h_{k}^{2}\sup_{t\in \lbrack 0,T]}\Vert u_{k}(t)\Vert _{L^{2}(\Omega
)}^{4}.  \label{eq:I3}
\end{equation}
It follows from (\ref{eq:I11}), (\ref{eq:I2}), (\ref{eq:I3}) and
 (\ref{eq:contradiction2}) that 
\begin{equation}
\lim_{k\rightarrow \infty }I_{k}=0.  \label{eq:limit}
\end{equation}
As the linear combination of $\chi(t,\tau)b(x,\xi)$ are dense in $\Con_0^\infty(T^\star(\RR^{d+1}))$, using (\ref{eq:carac}) and (\ref{eq:limit}), 
we deduce that $m_{0}=(x_{0},t_{0},\xi _{0},\tau _{0})\notin \supp \mu $.
\end{case}

\begin{case}
Assume that $x_{0}\in \partial \Omega .$

We would like to show that one can find a neighborhood $U_{x_{0}}$ of $x_{0}$
in ${\mathbb{R}}^{d}$ such that for any $b\in \Con_{0}^{\infty
}(U_{x_{0}}\times \mathbb{R}_{t}\times \mathbb{R}_{\xi }^{d}\times 
\mathbb{R}_{\tau }),$ we have 
\begin{equation}
\left\langle \mu ,(\tau +p)b\right\rangle =0 .  \label{eq:caraBord}
\end{equation}
Indeed this will imply that the point $m_{0}(x_{0},t_{0},\xi _{0},\tau _{0})$
(with $\tau _{0}+(x_{0},\xi _{0})\neq 0)$ does not belong to the support of 
$\mu $ as claimed.
Formula (\ref{eq:caraBord}) will be implied, by
\begin{equation}
\left\{ 
\begin{array}{l}
\lim\limits_{k\rightarrow +\infty }I_{k}=0\text{ where} \\ 
I_{k}=\left( b(x,t,h_{k}D_{x},h_{k}^{2}D_{t})\varphi
h_{k}^{2}(D_{t}+P)W_{k},W_{k}\right) _{L^{2}(\mathbb{R}^{d+1})}.
\end{array}
\right.  \label{eq:termeIk}
\end{equation}
where $\varphi \in \Con_{0}^{\infty }(U_{x_{0}}),\varphi =1$ on 
$\pi _{x}\supp b.$
Let $U_{x_{0}}$ a neighborhood of $x_0$ such that there exists a 
$\Con^{\infty }$\ diffeomorphisme $F$ from $U_{x_{0}}$\ to a neighborhood $U_{0}$
of the origin in $\mathbb{R}^{d}$ satisfying, 
\begin{equation}
\left\{ 
\begin{array}{c}
F(U_{x_{0}}\cap \Omega )=\{y\in U_{0}:y_{1}>0\} \\ 
F(U_{x_{0}}\cap \partial \Omega )=\{y\in U_{0}:y_{1}=0\} \\ 
(P(x,D)W_{k})\circ F^{-1}=(D_{1}^{2}+R(y,D^{\prime })
+ L(x,D))
(W_{k}\circ F^{-1}),
\end{array}
\right.  \label{eq:U0}
\end{equation}
where $R$ is a second-order differential operator, $D^{\prime
}=(D_{2},...,D_{d})$
and $L(x,D)$ a first order differential operator.
Let us set 
\begin{equation}
v_{k}=u_{k}\circ F^{-1},\text{ \ \ }V_{k}=1_{[0,T]}1_{y_{1}>0}v_{k},
\label{eq:v0}
\end{equation}
then we will have 
\begin{equation}
\left\{ 
\begin{array}{l}
\left( D_{t}+D_{1}^{2}+R(y,D^{\prime })+ L(x,D)\right) v_{k}=iaP^{1/2}(au_{k})\circ
F^{-1}+h_{k}^{-1}g_{k}\circ F^{-1}:=f_{k} \\ 
v_{k}|_{y_{1}=0}=0.
\end{array}
\right.  \label{eq:probRedresse}
\end{equation}
Making the change of variable $\ x=F^{-1}(y)$ on the right-hand side of the
second line of (\ref{eq:termeIk}), we see that 
\begin{equation*}
I_{k}=\left( \tilde{b}(y,t,h_{k}D_{y},h_{k}^{2}D_{t})\psi
h_{k}^{2}(D_{t}+D_{1}^{2}+R(y,D^{\prime })+ L(x,D))V_{k},V_{k}\right) _{L^{2}
(\mathbb{R}^{d+1})},
\end{equation*}
where $\tilde{b}\in \Con_{0}^{\infty }(U_{0}\times \mathbb{R}_{t}\times 
\mathbb{R}_{\eta }^{d}\times \mathbb{R}_{\tau }),$ and $\psi \in 
\Con_{0}^{\infty }(U_{0}),$ $\psi =1$ on $\pi _{y}\supp\tilde{b}.$
To prove (\ref{eq:termeIk}) it is sufficient to prove that,
\begin{equation*}
\lim_{k\rightarrow +\infty }J_{k}=\lim_{k\rightarrow +\infty }\left( T\psi
_{0}(y_{1})\psi _{1}(y^{\prime })h_{k}^{2}(D_{t}+D_{1}^{2}+R(y,D^{\prime
})+ L(x,D))V_{k},V_{k}\right) _{L^{2}(\mathbb{R}^{d+1})}=0,
\end{equation*}
where $T=\theta (y_{1},h_{k}D_{1})\Phi (y^{\prime },h_{k}D^{\prime })\chi
(t,h_{k}^{2}D_{t}),$ $\theta \Phi \chi \in \Con_{0}^{\infty }(U_{0}\times 
\mathbb{R}_{t}\times \mathbb{R}_{\eta }^{d}\times \mathbb{R}_{\tau }),$ 
$\psi _{0}\psi _{1}\in \Con_{0}^{\infty }(U_{0}),$ $\psi _{0}\psi _{1}=1$ on $
\pi _{y}\supp\theta \Phi \chi $;
According to (\ref{eq:probRedresse}) we have,
\begin{align*}
(D_{t}+D_{1}^{2}+R(y,D^{\prime })+ L(x,D))V_{k}& =f_{k}-i1_{y_{1}>0}v_{k}(0,.)\delta
_{t=0}+i1_{y_{1}>0}v_{k}(T,.)\delta _{t=T} \\
& \quad -i1_{[0,T]}(D_{1}v_{k}|_{y_{1}=0})\otimes \delta _{y_{1}=0}.
\end{align*}
Therefore (\ref{eq:termeIk}) will be proved if we can prove that
\begin{equation}
\left\{ 
\begin{array}{l}
\lim\limits_{k\rightarrow +\infty }A_{k}^{j}=0,\text{ \ }j=1,2,3,
\text{ where } \\ 
A_{k}^{1}=\left( \theta (y_{1},h_{k}D_{1})\Phi (y^{\prime },h_{k}D^{\prime
})\chi (t,h_{k}^{2}D_{t})\psi _{0}\psi
_{1}h_{k}^{2}1_{y_{1}>0}v_{k}(s,.)\delta _{t=s},V_{k}\right) ,\text{ }s=0,
\text{ }T, \\[3pt] 
A_{k}^{2}=\left( \theta (y_{1},h_{k}D_{1})\Phi (y^{\prime },h_{k}D^{\prime
})\chi (t,h_{k}^{2}D_{t})\psi _{0}\psi
_{1}h_{k}^{2}1_{[0,T]}(D_{1}v_{k}|_{y_{1}=0})\otimes \delta
_{y_{1}=0},V_{k}\right) , \\[3pt] 
A_{k}^{3}=\left( \theta (y_{1},h_{k}D_{1})\Phi (y^{\prime },h_{k}D^{\prime
})\chi (t,h_{k}^{2}D_{t})\psi _{0}\psi _{1}h_{k}^{2}f_{k},V_{k}\right) .
\end{array}
\right.  \label{eq:A-k}
\end{equation}
As in \cite[A.18]{RZ} 
\begin{equation}
\lim\limits_{k\rightarrow +\infty }A_{k}^{1}=0 .  \label{A1}
\end{equation}
To estimate the term $A_{k}^{2}$ we need a Lemma.
With $U_{0}$ introduced in (\ref{eq:U0}), we set $U_{0}^{+}=\{y\in
U_{0}:y_{1}>0\}.$ We consider a smooth solution of the problem:
\begin{equation}
\left\{ 
\begin{array}{l}
\left( D_{t}+D_{1}^{2}+R(y,D^{\prime })+ L(x,D)\right) u=g\text{ \ in \ }
U_{0}^{+}\times \mathbb{R}_{t} \\[2pt] 
u|_{y_{1}=0}=0
\end{array}
\right.  \label{eq:sec}
\end{equation}

\begin{lemma}
\label{LemmaA6} Let $\chi \in \Con_{0}^{\infty }(U_{0})$ and $\chi _{1}\in 
\Con_{0}^{\infty }(U_{0})$ $\chi _{1}=1$ on $\supp \chi .$ There exists $C>0$
such that for any solution $u$ of (\ref{eq:sec}) and all $h$ in $]0,1],$ we
have
\begin{align*}
\int_{0}^{T}\left\Vert \left( \chi h\partial _{1}u\right)
_{|y_{1}=0}(t)\right\Vert _{L^{2}}^{2}dt &\leq C\left(
\int_{0}^{T}\sum_{|\alpha |\leq 1}\left\Vert \chi _{1}(hD)^{\alpha
}u(t)\right\Vert _{L^{2}(U_{0}^{+})}^{2}dt\right. \\
&\quad +\left\Vert h^{\frac{1}{2}}\chi u(0)\right\Vert
_{L^{2}(U_{0}^{+})}\left\Vert h^{\frac{1}{2}}(h\partial _{1}u)(0)\right\Vert
_{L^{2}(U_{0}^{+})} \\
&\quad \left. +\left\Vert h^{\frac{1}{2}}\chi u(T)\right\Vert
_{L^{2}(U_{0}^{+})}\left\Vert h^{\frac{1}{2}}(h\partial _{1}u)(T)\right\Vert
_{L^{2}(U_{0}^{+})}+\left\Vert \chi _{1}hg\right\Vert _{L^{2}}^{2}\right).
\end{align*}
\end{lemma}
\end{case}
\begin{prooff}[Proof of the Lemma]
It is analogue to the proof of \cite[Lemma A.6]{RZ}.
\end{prooff}
We replace in the previous Lemma $g$ by $iaP^{1/2}(au_k)\circ
F^{-1}+h_k^{-1}g_k\circ F^{-1} $ and by \eqref{eq:v0}, we obtain easily the
following corollary.

\begin{corollary}
\label{corollaryA7} One can find a constant $C>0$ such that 
\begin{align*}
\int_{0}^{T}\left\Vert \left( \chi h_{k}\partial _{1}v_{k}\right)
_{|y_{1}=0}(t)\right\Vert _{L^{2}}^{2}dt &\leq C\left(
\int_{0}^{T}\left\Vert \tilde{\chi}u_{k}(t)\right\Vert _{L^{2}(\Omega
)}^{2}dt+\left\Vert h_{k}^{1/2}u_{k}(0)\right\Vert _{L^{2}(\Omega
)}^{2}dt\right. \\
&\quad \left. +\int_{0}^{T}\left( \left\Vert \tilde{\chi}
a(h_{k}^{2}P)^{1/2}au_{k}\right\Vert _{L^{2}}^{2}+\left\Vert \tilde{\chi}
g_{k}\right\Vert _{L^{2}}^{2}\right) dt\right) \\
&\leq C,
\end{align*}
where $v_{k}$ has been defined in (\ref{eq:v0}) and $\tilde{\chi}\in 
\Con_{0}^{\infty }(\mathbb{R}^{d}).$
\end{corollary}

Let us go back to the estimate of $A_{k}^{2}$ defined in (\ref{eq:A-k}). We
have
\begin{equation*}
\left\vert A_{k}^{2}\right\vert ^{2}\leq Ch_{k}^{2}\left\Vert \theta
(y_{1},h_{k}D_{1})\delta _{y_{1}=0}\right\Vert _{L^{2}
(\mathbb{R})}^{2}\left\Vert \left( \psi _{2}V_{k}\right) \right\Vert _{L^{2}
(\mathbb{R}^{d+1})}^{2}\int_{0}^{T}\left\Vert \left( \psi _{1}h_{k}D_{1}v_{k}\right)
_{|y_{1}=0}(t)\right\Vert _{L^{2}(\mathbb{R}^{d-1})}^{2}dt.
\end{equation*}
Applying (\ref{eq:bound}), \cite[Lemma A.5]{RZ} with $p=2$, $l=1$ and
corollary \ref{corollaryA7}, we obtain
\begin{equation}
\left\vert A_{k}^{2}\right\vert \leq ch_{k}\longrightarrow 0 .  \label{eq:A2}
\end{equation}
The term $\left\vert A_{k}^{3}\right\vert $ can be treated as the first and
the second term in the case 1.\newline
Using (\ref{A1}) and (\ref{eq:A2}), we deduce (\ref{eq:A-k}), which implies 
(\ref{eq:termeIk}) thus (\ref{eq:caraBord}). The proof of 
Proposition~\ref{prop:sopport}\ is complete. 
\end{prooff}

\subsection{The microlocal defect measure does not vanish identically}

First let us prove that the sequence $(u_k)$ have mass in a compact domain.

\begin{lemma}
\label{Lemma2.6} There exists a subsequence $k_\nu$, there exists $R>0$ such
that 
\begin{equation*}
\int_0^T\| u_{k_\nu}(t)\|_{L^2(x\in\Omega,\ |x|<R)}^2dt\ge 1/2.
\end{equation*}
\end{lemma}
\begin{prooff}[Proof of Lemma]
We prove the Lemma by contradiction. Assume that 
\begin{equation}
\forall R>R_{0},\ \limsup_{k}\int_{0}^{T}\Vert u_{k}(t)\Vert _{L^{2}(x\in
\Omega ,\ |x|\leq 2R
+1)}^{2}dt\leq 3/4,  \label{Contradiction lemme 2.6}
\end{equation}
where $R_{0}$ is large enough such that $\supp a\subset \{|x|\leq R_{0}/2\}$.

Let $\chi \in \Con^{\infty }(\mathbb{R}^{d})$ such that $\chi =1$ for $|x|>2$
and $\chi =0$ for $|x|<1$. We set $\chi _{R}(x)=\chi (x/R)$ and by the
choice of $R_{0}$ we have $a\chi _{R}=\chi _{R}a=0$ . The function
$v_{k}:=\chi _{R}u_{k}$ satisfies 
\begin{equation*}
D_{t}v_{k}+Pv_{k}=h_{k}^{-1}\chi _{R}g_{k}+[P,\chi _{R}]u_{k}.
\end{equation*}
From \cite[Theorem 2.8]{doi05}, we have 
\begin{equation}
\int_{0}^{T}\Vert \langle x\rangle ^{-s}v_{k}\Vert _{L^{2}
(\mathbb{R}^{d})}^{2}\leq C(\Vert E_{-\frac{1}{2}}v_{k}(0)\Vert _{L^{2}
(\mathbb{R}^{d})}^{2}+\int_{0}^{T}\Vert \langle x\rangle ^{s}E_{-1}\left(
h_{k}^{-1}\chi _{R}g_{k}+[P,\chi _{R}]u_{k}\right) \Vert _{L^{2}
(\mathbb{R}^{d})}^{2}dt),  \label{Doi}
\end{equation}
where $E_{s}$ is the pseudo-differential operator with symbol $e_{s}=(1+{p}
(x,\xi )+|x|^{2})^{\frac{s}{2}}$ which belongs to $S((|\xi |+<x>)^{s},g)$. 

For the first term of the right hand side of (\ref{Doi}) we have, 
 where $(\cdot,\cdot)$ means the scalar product in $L^2(\Omega)$,

\begin{align*} 
\Vert E_{-\frac{1}{2}}v_{k}(0)\Vert _{L^{2}}^{2} &=h_k\Vert E_{-\frac{1}{2}
}\chi_R P^{\frac{1}{4}}(h_k^{2}P)^{-\frac{1}{4}}\psi _{1}(h_k^{2}P)\psi
(h_k^{2}P)u(0)\Vert _{L^{2}}^{2}, \\
&=h_k ( S\psi _{2}(h_k^{2}P)u_{k}(0),S\psi
_{2}(h_k^{2}P)u_{k}(0)) ,\text{ where } S=E_{-\frac{1}{2}}\chi_R P^{
\frac{1}{4}},\text{ and }\psi _{2}(t)=t^{-\frac{1}{4}}\psi _{1} \\
&=h_k( \psi _{2}(h_k^{2}P)S^{\star }S\psi
_{2}(h_k^{2}P)u_{k}(0),u_{k}(0))  \\
&=h_k( \psi _{2}(h_k^{2}P)(h_k^{2}P)^{-\frac{1}{4}}Q\chi_R (h_k^{2}P)^{
\frac{1}{4}}\psi _{2}(h_k^{2}P)u_{k}(0),u_{k}(0)) \\
&\leq Ch_k\Vert u_{k}(0)\Vert _{L^{2}}^{2},
\end{align*}
where $\psi _{1}\in \Con_{0}^{\infty }(0,+\infty) $ and $\psi _{1}=1\text{
on }\supp(\psi )$, $S^{\star }S=P^{-\frac{1}{4}}Q\chi_R P^{\frac{1}{4}}$, $
Q=P^{\frac{1}{2}}\chi_R A_{-1}$, and $A_{-1}=E_{-\frac{1}{2}}^{\star }
E_{-\frac{1}{2}}$. We have used that the operator $Q$ is bounded from $L^{2}(
\mathbb{R}^{d})$ to $L^{2}(\Omega )$ (see \cite[Lemma 4.2]{RZ}).\newline
Then from (\ref{eq:contradiction2}), we deduce that 
\begin{equation}
\lim_{k\rightarrow +\infty }\Vert E_{-\frac{1}{2}}v_{k}(0)\Vert
_{L^{2}}^{2}=0 .  \label{eq:E-1/2}
\end{equation}
Concerning the term $\displaystyle\int_{0}^{T}\Vert \langle x\rangle
^{s}E_{-1}h_k^{-1}\chi_R 
g_{k}\Vert _{L^{2}}^{2}dt$, we will prove that it
tends to zero.\newline
Let $\psi _{1}\in \Con_{0}^{\infty }(\mathbb{R})$, such that $\psi _{1}=1$
on $\supp\psi $.\newline
Since $\psi _{1}(h_{k}^{2}P)u_{k}=u_{k}$ then applying $1-\psi_1(h_k^2P)$ to
Formula \eqref{eq:lowfre2}, we obtain 
\begin{equation*}
h_{k}^{-1}g_{k}=h_{k}^{-1}\psi
_{1}(h_{k}^{2}P)g_{k}-ih_{k}^{-1}a(h_k^{2}P)^{1/2}a\psi
_{1}(h_{k}^{2}P)u_{k}+ih_{k}^{-1}\psi
_{1}(h_{k}^{2}P)a(h_k^{2}P)^{1/2}au_{k}.
\end{equation*}
Using that $\chi_R a=0$, we have 
\begin{equation*}
h_{k}^{-1}\chi_R g_{k}=h_{k}^{-1}\chi_R \psi
_{1}(h_{k}^{2}P)g_{k}+ih_{k}^{-1}\chi_R \psi
_{1}(h_{k}^{2}P)a(h_k^{2}P)^{1/2}au_{k}.
\end{equation*}
And then 
\begin{align*}
&\int_{0}^{T}\Vert \langle x\rangle ^{s}E_{-1}h_{k}^{-1}\chi_R g_{k}\Vert
_{L^{2}}^{2}dt \\
&\quad \quad \leq \int_{0}^{T}\Vert \langle x\rangle ^{s}E_{-1}\chi_R
h_k^{-1}\psi _{1}(h_{k}^{2}P)g_{k}\Vert ^{2}dt+\int_{0}^{T}\Vert \langle
x\rangle ^{s}E_{-1}\chi_R h_{k}^{-1}\psi
_{1}(h_{k}^{2}P)a(h_{k}^{2}P)^{1/2}au_{k}\Vert ^{2}dt \\
&\quad \quad \leq \int_{0}^{T}\Vert \langle x\rangle ^{s}E_{-1}\chi_R
P^{1/2}\psi _{2}(h_{k}^{2}P)g_{k}\Vert ^{2}dt+\int_{0}^{T}\Vert \langle
x\rangle ^{s}E_{-1}\chi_R h_{k}^{-1}\psi
_{1}(h_{k}^{2}P)a(h_{k}^{2}P)^{1/2}au_{k}\Vert ^{2}dt,
\end{align*}
where $\psi _{2}(t)=t^{-1/2}\psi _{1}(t).$
We have, 
\begin{align}
\int_{0}^{T}\Vert \langle x\rangle ^{s}E_{-1}\chi_R P^{1/2}\psi
_{2}(h_{k}^{2}P)g_{k}\Vert ^{2}dt&\le I + I\!I, 
\end{align}
where 
\begin{equation}
I= \int_{0}^{T}\Vert \langle x\rangle^{s}E_{-1}\langle x\rangle ^{-s}\chi_R P^{1/2}\psi _{2}(h_{k}^{2}P)\langle x\rangle ^{s}g_{k}\Vert ^{2}dt  \notag
\end{equation}
and 
\begin{equation}
I\!I=h_k^{-2} \int_{0}^{T}\Vert \langle x\rangle ^{s}E_{-1}\chi_R
[(h_k^2P)^{1/2}\psi _{2}(h_{k}^{2}P),\langle x\rangle ^{-s} ]\langle
x\rangle ^{s}g_{k}\Vert ^{2}dt.  \notag 
\end{equation}
It follows that the symbol of $\langle x\rangle ^{s}E_{-1}\langle x\rangle ^{-s}$
belongs to $S((|\xi|+\langle x\rangle)^{-1})$ then $\langle x\rangle
^{s}E_{-1}\langle x\rangle ^{-s}\chi_R 
P^{1/2}$ is bounded on $L^2(\Omega)$
(see \cite[Lemma 4.2]{RZ}) and we have
\begin{equation}
I\le C \int_{0}^{T}\Vert \langle x\rangle ^{s}g_{k}\Vert ^{2}dt ,  \notag
\end{equation}
According to Lemma~\ref{lemma:D}, $h_k^{-1}\langle x\rangle ^{s}
[(h_k^2P)^{1/2}\psi_{2}(h_{k}^{2}P),\langle x\rangle ^{-s} ]$ is bounded on 
$L^2(\Omega)$ and we get
\begin{equation}
I\!I\le C \int_{0}^{T}\Vert \langle x\rangle ^{s}g_{k}\Vert ^{2}dt.  \notag
\end{equation}
To estimate $$\int_{0}^{T}\Vert \langle x\rangle ^{s}E_{-1}\chi
_{R}h_{k}^{-1}\psi _{1}(h_{k}^{2}P)a(h_{k}^{2}P)^{1/2}au_{k}\Vert ^{2}dt,$$
we have with $\psi _{2}(s)=s^{-1}\psi
_{1}(s)$ and $\tilde{\chi}$ a smooth function such that, $\tilde{\chi}=1$
for $|x|\geq 1$ and $\tilde{\chi}=0$ for $|x|\leq 1/2$, $\tilde{\chi}_{R}(x)=
\tilde{\chi}(x/R)$, 
\begin{align}
\langle x\rangle ^{s}E_{-1}\chi _{R}h_{k}^{-1}\psi _{1}(h_{k}^{2}P)a&
=\langle x\rangle ^{s}E_{-1}\chi _{R}Ph_{k}\psi _{2}(h_{k}^{2}P)a=\langle
x\rangle ^{s}E_{-1}\chi _{R}P\tilde{\chi}_{R}h_{k}\psi _{2}(h_{k}^{2}P)a 
\notag \\
& =\langle x\rangle ^{s}E_{-1}\langle x\rangle ^{-s}\chi
_{R}P^{1/2}(h_{k}^{2}P)^{1/2}\langle x\rangle ^{s}[\tilde{\chi}_{R},\psi
_{2}(h_{k}^{2}P)]a  \notag \\
& \quad +\langle x\rangle ^{s}E_{-1}\langle x\rangle ^{-s}\chi _{R}[\langle
x\rangle ^{s},P]\tilde{\chi}_{R}h_{k}[\psi _{2}(h_{k}^{2}P),a],
\label{egalite pour le terme 2}
\end{align}
where we have used $a\tilde{\chi}_{R}=0$ if $R$ large enough.

By the \cite[Lemma A.5]{RZ} and Lemma~\ref{lemma:B bis} the first term of 
\eqref{egalite pour le terme 2} is bounded on $L^2(\Omega)$ by $Ch_k$. As 
$[\langle x\rangle ^{s},P]$ is a sum of term $\alpha \partial_{x_j}$ where 
$\alpha $ is bounded, $\langle x\rangle ^{s}E_{-1}\langle x\rangle
^{-s}\chi_R[\langle x\rangle ^{s},P]$ is bounded on $L^2(\Omega)$, and 
$[\psi_{2}(h_{k}^{2}P),a]$ is bounded on $L^2(\Omega)$ by \cite[Lemma 6.3]{RZ}. Then the second term of \eqref{egalite pour le terme 2} is bounded on 
$L^2(\Omega)$ by $Ch_k$. Finally, we yield 
by Lemma~\ref{lemma:H},

\begin{align}
\int_{0}^{T}\Vert \langle x\rangle ^{s}E_{-1}\chi_R h_{k}^{-1}\psi
_{1}(h_{k}^{2}P)a(h_{k}^{2}P)^{1/2}au_{k}\Vert ^{2}dt &\leq C_R
h_{k}^{2}\int_{0}^{T}\Vert 
(h_{k}^{2}P)^{1/2}a
 u_{k}\Vert ^{2}dt  \notag \\
&
\le C_Rh_k^2\sup_{t\in [0,T]}\|u_k(t,.)\|^2.
\end{align}
According to (\ref{eq:contradiction}) and (\ref{eq:contradiction2}), we
conclude that the second term of the right hand side of (\ref{Doi}) goes to
zeros when $k$ tend to $+\infty $
\begin{equation}
\lim_{k\rightarrow \infty }\int_{0}^{T}\Vert \langle x\rangle
^{s}E_{-1}h_{k}^{-1}\chi _{R}g_{k}\Vert _{L^{2}}^{2}dt=0.
\label{Doi terme 2}
\end{equation}
Now we estimate the term $\displaystyle\int_{0}^{T}\Vert \langle x\rangle
^{s}E_{-1}[P,\chi _{R}]u_{k}\Vert _{L^{2}}^{2}dt$.\newline
Let $\chi _{1}\in \Con_{0}^{\infty }(R-1<|x|<2R+1),\chi _{1}\geq 0,\chi
_{1}=1\mbox{ on
}\supp(\nabla \chi _{R}),$ 
\begin{align}
\int_{0}^{T}\Vert \langle x\rangle ^{s}E_{-1}[P,\chi _{R}]u_{k}\Vert
_{L^{2}}^{2})dt& \leq \int_{0}^{T}\Vert \langle x\rangle ^{s}\chi
_{1}E_{-1}[P,\chi _{R}]\chi _{1}u_{k}]\Vert _{L^{2}(\Omega )}^{2}dt  \notag
\\
& \quad +\int_{0}^{T}\Vert \langle x\rangle ^{s}(1-\chi _{1})E_{-1}[P,\chi
_{R}]\chi _{1}u_{k}]\Vert _{L^{2}(\Omega )}^{2}dt,\qquad  \notag \\
& \leq CR^{2(s-1)}\int_{0}^{T}\Vert u_{k}\Vert
_{L^{2}(R-1<|x|<2R+1)}^{2}dt\leq CR^{2(s-1)},  \label{Doi terme 3}
\end{align}
where we have used, first that $E_{-1}\partial _{x}$ is bounded on $L^{2}$, $
\langle x\rangle ^{s}$ is estimate by $CR^{s}$ on support of $\chi _{1}$ and 
$\partial _{x}\chi _{R}$ is the product of a bounded function by $R^{-1}$,
second, the symbol of $\langle x\rangle ^{s}(1-\chi _{1})E_{-1}[P,\chi _{R}]$
is uniformly bounded in $R^{-1}S((\langle x\rangle +|\xi |)^{-N},g)$ for all 
$N$. The last inequality uses the contradiction assumption 
\eqref{Contradiction lemme 2.6}.

Following \eqref{Doi}, \eqref{eq:E-1/2}, \eqref{Doi terme 2} and 
\eqref{Doi
terme 3}, we have, 
\begin{equation*}
\int_{0}^{T}\Vert \langle x\rangle ^{-s}u_{k}\Vert
_{L^{2}(|x|>2R)}^{2}dt\leq \int_{0}^{T}\Vert \langle x\rangle
^{-s}v_{k}\Vert _{L^{2}(\mathbb{R}^{d})}^{2}\le C_R\delta_k+CR^{2(s-1)},
\end{equation*}
where $\delta_k\to 0$ when $k\to +\infty $, $C$ is independent of $R$ and $
C_R$ may depend of $R$. Then we have 
\begin{align*}
\int_0^T\| u_k\|^2_{L^2(x\in\Omega,\ |x|<2R)}&\ge \int_0^T\| \langle
x\rangle^{-s}u_k\|^2_{L^2(x\in\Omega,\ |x|<2R)} \\
&\ge \int_0^T\| \langle x\rangle^{-s}u_k\|^2_{L^2(x\in\Omega)} -\int_0^T\|
\langle x\rangle^{-s}u_k\|^2_{L^2( |x|>2R)} \\
&\ge \int_0^T\| \langle
x\rangle^{-s}u_k\|^2_{L^2(x\in\Omega)}-C_R\delta_k-CR^{2(s-1)}.
\end{align*}
This with \eqref{eq:contradiction3} implies a contradiction with 
\eqref{Contradiction lemme 2.6} and proves the Lemma.
\end{prooff}

In the sequel, for simplicity, we
shall denote
 the sequence $u_{k_\nu}$ found in
Lemma~\ref{Lemma2.6} by $u_k$. Thus there exist $R_{0}>0$, $k_{0}>0$ such
that

\begin{equation*}
\int_{0}^{T}\Vert u_{k}(t)\Vert _{L^{2}(|x|<R)}^{2}dt\geq \frac{1}{2},
\end{equation*}
when $R>R_{0}$ and $k>k_{0}$. \newline
We consider $\chi _{1}\in \Con_{0}^{\infty }(\mathds{R}^{d})$ such that $$
0\leq \chi _{1}\leq 1,\; \chi _{1}(x)=1 \mbox{ if }|x|\leq R_{1}+2\mbox{ and }\supp\chi
_{1}\subset \{|x|\leq R_{1}+3\},$$ with $R_{1}>R_{0}$.\\ Let $A\geq 1$, $R\geq
1 $, $\psi _{A}\in \Con_{0}^{\infty }(\mathds{R})$, $\phi _{R}\in 
\Con_{0}^{\infty }(\mathds{R})$ be such that $0\leq \psi _{A}$, 
$\phi _{R}\leq 1$
and 
\begin{equation*}
\psi _{A}(\tau )=1\text{ if }|\tau |\leq A,\phi _{R}(t)=1\text{ if }|t|\leq
R.   
\end{equation*}
We recall that $w_k(t)=1_\Omega u_k(t)$.

\begin{proposition}
There exist positive constants $A_0$, $R_0$, $k_0$ such that 
\begin{equation*}
\int_{\mathds{R}}\|\psi_A(h^2_kD_t)\phi_R(h^2_k\Delta)1_{[0,T]}
\chi_1w_k(t)\|^2_{L^2(\mathds{R}^d)}dt\geq \frac{1}{4},
\end{equation*}
when $A\geq A_0$, $R\geq R_0$, $k\geq k_0$.
\end{proposition}

\begin{corollary}
The measure $\mu$ does not vanish identically.
\end{corollary}
\begin{prooff}[Proof of Proposition ]
Set $I=(Id-\psi _{A}(h_{k}^{2}D_{t}))1_{[0,T]}\chi _{1}u_{k}$ and 
$\widetilde{\psi }(\tau )=\dfrac{1-\psi _{A}(\tau )}{\tau }$. It is easy to
see that $\widetilde{\psi }\in L^{\infty }(\mathbb{R})$ and 
$|\widetilde{\psi }(\tau )|\leq \frac{1}{A}$ for all 
$\tau \in \mathbb{R}$.

We have 
\begin{align*}
I &=\widetilde{\psi }_{A}(h_{k}^{2}D_{t})h_{k}^{2}D_{t}(1_{[0,T]}\chi
_{1}w_{k}) \\
&=\frac{h_k^{2}}{i}\widetilde{\psi }_{A}(h_{k}^{2}D_{t})\chi
_{1}(u_{k}(0)\delta _{t=0}-u_{k}(T)\delta _{t=T}) \\
&\quad \widetilde{\psi }_{A}(h_{k}^{2}D_{t})\chi
_{1}1_{[0,T]}(-h_{k}^{2}Pu_{k}+ih_{k}a(h_{k}^{2}P)^{1/2}
au_{k}+h_{k}g_{k}) \\
&=B_{k}^{1}+B_{k}^{2}+B_{k}^{3}+B_{k}^{4}.
\end{align*}
From \cite[See the proof of Proposition 6.1]{RZ} we know that $\Vert 
\widetilde{\psi }_{A}(h_{k}^{2}D_{t})\delta _{t=a}\Vert _{L^{2}
(\mathbb{R})}\leq Ch_{k}^{-1}$, then we deduce that 
\begin{equation*}
\lim_{k\rightarrow +\infty }\int_{\mathbb{R}}\Vert B_{k}^{1}\Vert
_{L^{2}(\Omega )}^{2}dt\leq \lim_{k\rightarrow +\infty
}Ch_{k}^{4}h_{k}^{-2}(\Vert u_{k}(0)\Vert _{L^{2}(\Omega )}^{2}+\Vert
u_{k}(T)\Vert _{L^{2}(\Omega )}^{2})=0.
\end{equation*}
Using (\ref{eq:auk}) and (\ref{eq:contradiction2}), we can prove easily that 
\begin{equation*}
\lim_{k\rightarrow +\infty }\int_{\mathbb{R}}\Vert B_{k}^{3}\Vert
_{L^{2}(\Omega )}^{2}dt\leq C\lim_{k\rightarrow +\infty
}\int_{0}^{T}h_{k}\Vert (h_{k}^{2}P)^{1/2}au_{k}\Vert _{L^{2}(\Omega
)}^{2}dt=0 .  
\end{equation*}
From \eqref{eq:contradiction} we can see that 
\begin{equation*}
\lim_{k\rightarrow +\infty }\int_{\mathbb{R}}\Vert B_{k}^{4}\Vert
_{L^{2}(\Omega )}^{2}dt\leq C\lim_{k\rightarrow +\infty }\int_{0}^{T}\Vert
\chi _{1}g_{k}\Vert _{L^{2}(\Omega )}^{2}dt=0 .  
\end{equation*}
Now, for $B_{k}^{2}$ we argue as in \cite[See the proof of 
Proposition 6.1]{RZ}. Let $\tilde\theta\in \Con_0^\infty (0,+\infty) $ 
such $\tilde\theta=1$ on the support of $\psi$ and let 
$\tilde\theta_1(s)=s\tilde\theta(s)$. We
have 
\begin{align*}
B_{k}^{2} &=-\widetilde{\psi }_{A}(h_{k}^{2}D_{t})
\chi _{1}1_{[0,T]}h_k^{2}P\widetilde{\theta }(h_{k}^{2}P)u_{k} \\
&=-\widetilde{\psi }_{A}(h_{k}^{2}D_{t})1_{[0,T]}[\chi _{1},
\widetilde{\theta }_{1}(h_{k}^{2}P)]u_{k}-\widetilde{\psi }_{A}
(h_{k}^{2}D_{t})1_{[0,T]}
\widetilde{\theta }_{1}(h_k^{2}P)\chi _{1}u_{k}.
\end{align*}
Using Lemma 6.3 in \cite{RZ} and the fact that 
\begin{equation*}
\Vert \widetilde{\psi }_{A}(h_{k}^{2}D_{t})\Vert _{L^{2}(\mathbb{R}
)\rightarrow L^{2}(\mathbb{R})}=O\left( \frac{1}{A}\right) ,\,\Vert 
\widetilde{\theta }_{1}(h_{k}^{2}P)\Vert _{L^{2}
(\Omega)\rightarrow
L^{2}(\Omega})=O(1),
\end{equation*}
uniformly in $k$, we deduce that 
\begin{equation*}
\int_{\mathbb{R}}\Vert B_{k}^{2}\Vert _{L^{2}(\Omega )}^{2}dt\leq
C(h_{k}^{2}\sup_{t\in \lbrack 0,T]}\Vert u_{k}(t)\Vert _{L^{2}(\Omega
)}^{2}dt+\frac{1}{A}\int_{0}^{T}\Vert \chi _{1}u_{k}\Vert _{L^{2}(\Omega
)}^{2}dt). 
\end{equation*}
Taking $k$ and $A$ sufficiently large we obtain 
\begin{equation}
\int_{\mathds{R}}\Vert \psi _{A}(h_{k}^{2}D_{t})1_{[0,T]}\chi
_{1}w_{k}(t)\Vert _{L^{2}(\mathds{R}^{d})}^{2}dt\geq \frac{1}{3}.
\label{trocature D_t}
\end{equation}
Now, we set 
\begin{equation*}
\text{II} =(Id-\phi _{R}(h_{k}^{2}\Delta ))\psi
_{A}(h_{k}^{2}D_{t})1_{[0,T]}\chi _{1}w_{k}.
\end{equation*}
It is proved in \cite{RZ} that 
\begin{equation}
\int_{\mathbb{R}}\Vert \text{II}\Vert _{L^{2}(\mathds{R}^{d})}^{2}dt\leq 
\frac{C_{R_1}}{R}(1+h_{k}^{2}),  \label{RZ1}
\end{equation}
where $C_{R_1}$ depends on $R_1$ and The proof does not depend on the
equation, so it remains valid in our case. Nevertheless we recall the proof
in the sequel for the convenience of the reader. Before we give the end of
the proof of proposition 2.7.

Taking $R$ sufficiently large and using (\ref{trocature D_t}), we obtain 
\begin{equation*}
\int_{\mathds{R}}\Vert \phi _{R}(h_{k}^{2}\Delta )\psi
_{A}(h_{k}^{2}D_{t})1_{[0,T]}\chi _{1}w_{k}(t)\Vert _{L^{2}(\mathds{R}
^{d})}^{2}dt\geq \frac{1}{4}.  
\end{equation*}

Return to the proof of \eqref{RZ1}. We have 
$\displaystyle |1-\phi_R(t)|\le C\frac{h|\xi|}{\sqrt{R}} $ then we obtain, 

\begin{align}
\int_{\mathbb{R}}\Vert \text{II}\Vert _{L^{2}(\mathds{R}^{d})}^{2}dt &\leq 
C\frac{h_{k}^{2}}{R}\int_{\mathbb{R}}\sum_j\Vert \partial _{j}\psi
_{A}(h_{k}^{2}D_{t})1_{[0,T]}\chi _{1}w_{k}\Vert _{L^{2}(\mathds{R}
^{d})}^{2}dt  \notag \\
&\leq C\frac{h_{k}^{2}}{R}\int_{\mathbb{R}}\sum_j\Vert \partial _{j}\psi
_{A}(h_{k}^{2}D_{t})1_{[0,T]}\chi _{1}u_{k}\Vert _{L^{2}(\Omega )}^{2}dt 
\notag \\
&\le \frac{h_{k}^{2}}{R}\sum_j\left( \int_{\mathbb{R}}\Vert \partial _{j}
\widetilde{\theta }(h_{k}^{2}P)\psi _{A}(h_{k}^{2}D_{t})1_{[0,T]}\chi
_{1}u_{k}\Vert _{L^{2}(\Omega )}^{2}dt\right.  \notag \\
&\quad +\left. \int_{\mathbb{R}}\Vert \partial _{j}(1-\widetilde{\theta }
(h_{k}^{2}P))\psi _{A}(h_{k}^{2}D_{t})1_{[0,T]}\chi _{1}u_{k}\Vert
_{L^{2}(\Omega )}^{2}dt\right)  \notag \\
&:=\frac{h_{k}^{2}}{R}(C_{k}^{1}+C_{k}^{2}),  \label{C^1_k+C^2_k}
\end{align}
where $\widetilde{\theta }\in \Con_{0}^{\infty }(\mathbb{R})$ satisfying
 $\widetilde{\theta }(t)=1$ if $t\in \supp(\theta _{1})$ and 
$\widetilde{\theta }\theta _{1}=\theta _{1}$.

We have by Lemma 6.3 \cite{RZ}
\begin{equation}
C^1_k \leq Ch_k^{-2}\int_0^T\| \chi_1 u_k\|^2_{L^2(\Omega)}dt\leq ch_k^{-2} ,
\label{C^1_k}
\end{equation}
and 
\begin{align}
C^2_k &\leq \int_{\mathbb{R}} \|\partial_j[\widetilde{\theta}
(h_k^2P),\chi_1]\psi_A(h_k^2D_t)1_{[0,T]}\widetilde{\chi}_1
u_k\|^2_{L^2(\Omega)}dt  \notag \\
&\leq \int_{\mathbb{R}} \|\psi_A(h_k^2D_t)1_{[0,T]}\widetilde{\chi}_1
u_k\|^2_{L^2(\Omega)}dt  \notag \\
&\leq C\int_0^T\| \widetilde{\chi}_1 u_k\|^2_{L^2(\Omega)}dt\leq C_{R_1}
\int_0^T\|\langle x\rangle^{-s} u_k\|^2_{L^2(\Omega)}dt,  \label{C^2_k}
\end{align}
where $\widetilde{\chi}_1\in \Con^\infty_0(\overline{\Omega})$, $\widetilde{
\chi}_1=1$ on $\supp (\chi_1)$.\newline
Combining (\ref{C^1_k+C^2_k}), (\ref{C^1_k}) and (\ref{C^2_k}), 
we obtain (\ref{RZ1}).
\end{prooff}

\subsection{The microlocal defect measure vanishes in the incoming set}

In this section we prove that the microlocal defect measure $\mu$ vanishes
in the incoming set.

First remind some notation introduced in \cite{RZ} section 7. We keep the
same notation when it is possible.

We denote by 
\begin{equation*}
b(x,\xi )=\sum_{j,k=1}^{d}b^{jk}(x)x_{j}\xi _{k}.
\end{equation*}

\begin{proposition}
\label{Prop Incom} Let $m_0=(x_0,t_0,\xi_0,\tau_0)\in T^\star(\R^{d+1})$ be such 
$\xi_0\not=0$, $\tau_0+p(x_0,\xi_0)=0$, $|x_0|\ge3R_0$, $b(x_0,\xi_0)\le
-3\delta|x_0||\xi_0|$ for some $\delta>0$ small enough. Then 
$m_0\notin\supp \mu$.
\end{proposition}

We remind the results proved in \cite{RZ} in section 7, Lemma 7.5 and
Corollary 7.6. A part of the proof is in Doi~\cite{Do}. We use the Weyl
quantification of symbol which is denoted by $Op^{w}$.

There exist a symbol $\Phi \in S(1,g)$ such that $0\leq \Phi \leq 1$ and a
symbol $\lambda _{1}\in S(1,g)$ such that, 
\begin{gather}
\supp\lambda _{1}\subset \supp\Phi \subset \{(x,\xi )\in T^{\ast }(\R^{d}),\
|x|\geq 2R_{0},\ b(x,\xi )\leq -\frac{\delta }{2}|x||\xi |,\ |\xi |\geq 
\frac{|\xi _{0}|}{4}\} ,  \label{Incom eq1} \\
\{(x,\xi )\in T^{\ast }(\R^{d}),\ |x|\geq \frac{5}{2}R_{0},\ b(x,\xi )\leq
-\delta |x||\xi |,\ |\xi |\geq \frac{|\xi _{0}|}{2}\}\subset \{(x,\xi )\in
T^{\ast }(\R^{d}),\ \Phi (x,\xi )=1\} ,  \notag \\     
\Phi (x,h\xi )=\Phi (x,\xi )\text{ when }|h\xi |\geq \frac{|\xi _{0}|}{2},
\text{ and }0<h\leq 1,  \notag \\
H_{{p}}\Phi (x,\xi )\leq 0\text{ on the support of }\lambda _{1} ,  \notag \\
\lambda _{1}\geq 0 ,  \notag \\
[\tilde{P},Op^{w}(\lambda _{1})]-\frac{1}{i}Op^{w}(H_{{p}}\lambda _{1})\in
Op^w(S(1,g)) ,  \label{Icom eq7} \\
\text{there exist two positive constants }C,\ C^{\prime }\text{ such that} ,
\notag \\
-H_{{p}}\lambda _{1}\geq C\langle x\rangle ^{-2s}\Phi ^{2}(x,\xi )(|x|+|\xi
|)-C^{\prime }\Phi^2(x,\xi ).  \label{Icom eq8}
\end{gather}

\begin{prooff}
Let $\varphi _{1}\in \Con_{0}^{\infty }(\R^{d})$ such that 
\begin{equation}
\varphi _{1}(x)=1\text{ if }|x|\leq \frac{4}{3}R_{0},\ \supp\varphi
_{1}\subset \{x,\ |x|\leq \frac{3}{2}R_{0}\}.  \label{Incom eq9}
\end{equation}
Let $M$ large enough such that, 
\begin{equation*}
|((1-\varphi _{1})Op^w(\lambda _{1})(1-\varphi _{1})u|u)|\leq \frac{M}{2}
\Vert u\Vert ^{2}.
\end{equation*}
Here and in the sequel $(\cdot |\cdot )$ and $\Vert \cdot \Vert $ denote the 
$L^{2}(\Omega )$ inner product and norm respectively. The cutoff make sense
with this $L^{2}$ product. We set, 
\begin{equation*}
N(t)=((M-(1-\varphi _{1}){\mathcal{O}}p(\lambda _{1})(1-\varphi
_{1}))u_{k}(t)|u_{k}(t)),
\end{equation*}
and we have 
\begin{equation}
\frac{M}{2}\Vert u_{k}(t)\Vert ^{2}\leq N(t)\leq 2M\Vert u_{k}(t)\Vert ^{2}.
\label{Incom eq11}
\end{equation}
Setting $\Lambda =M-(1-\varphi_1){\mathcal{O}} p(\lambda_1)(1-\varphi_1)$,
we have, 
\begin{equation*}
\frac d{dt}N(t)=(\Lambda \frac d{dt}u_k(t)|u_k(t))+(\Lambda u_k(t)| \frac
d{dt}u_k(t)).
\end{equation*}
From \eqref{eq:lowfre2} we have 
\begin{equation*}
\frac{d}{dt}
u_{k}=-iPu_{k}-h_{k}^{-1}a(h_{k}^{2}P)^{1/2}(au_{k})+ih_{k}^{-1}g_{k}.
\end{equation*}
We obtain, 
\begin{align}
\frac{d}{dt}N=& (i[P,\Lambda ]u_{k}|u_{k})  \notag \\
& -h_{k}^{-1}(\Lambda a(h^{2}P)^{1/2}au_{k}|u_{k})-h_{k}^{-1}(\Lambda
u_{k}|a(h_{k}^{2}P)^{1/2}au_{k})  \notag \\
& +ih_{k}^{-1}(\Lambda g_{k}|u_{k})-ih_{k}^{-1}(\Lambda u_{k}|g_{k})  \notag
\\
=& A_{1}+A_{2}+A_{3} .  \label{Icom eq16}
\end{align}
For support reasons, we have $a(1-\varphi_1)=0$ thus we deduce,

\begin{align}
A_{2}& =-\frac{M}{h_{k}}
[(a(h_{k}^{2}P)^{1/2}(au_{k})|u_{k})+(u_{k}|a(h_{k}^{2}P)^{1/2}(au_{k}))] 
\notag \\
& =-\frac{2M}{h_{k}}\Vert (h_{k}^{2}P)^{1/4}(au_{k})\Vert ^{2}\leq 0.
\label{Icom eq17}
\end{align}
We have, for a constant $C_1>0$ 
\begin{align}  \label{Icom eq18}
|A_3|\le \frac{C_1}{h_k}\|\langle x\rangle ^sg_k\|\| \langle x\rangle
^{-s}u_k\|.
\end{align}
To estimate $A_{1}$ we remark that $[P,\Lambda ]=[\tilde{P},\Lambda ]$ and 
\begin{equation}
\lbrack P,\Lambda ]=[\tilde{P},\varphi _{1}]{O}p^w(\lambda _{1})(1-\varphi
_{1})-(1-\varphi _{1})[\tilde{P},{O}p^w(\lambda _{1})](1-\varphi
_{1})+(1-\varphi _{1}){O}p^w(\lambda _{1})[\tilde{P},\varphi _{1}].
\label{Icom eq19}
\end{equation} 
Following \eqref{Incom eq1} and \eqref{Incom eq9}, the support of $\lambda_1$
and $\varphi_1$ are disjoint, thus, taking account of \eqref{Incom eq11}, we
have

\begin{equation}  \label{Incom eq20}
|(\big[[\tilde P,\varphi_1]{O} p^w(\lambda_1)(1-\varphi_1)+(1-\varphi_1){O}
p^w(\lambda_1)[\tilde P,\varphi_1]\big]u_k|u_k)|\le C_2N(t).
\end{equation}  
Let $d(x,\xi)\in \Con_{0}^{\infty }(\R^{2d}) $ supported in $\{|x-x_0|\le 1,\ |\xi-\xi_0|\le 1\}$, and 
$d(x_0,\xi_0)=1$. According to \eqref{Icom eq7}, \eqref{Icom eq8} and G\aa rding
inequality, we get,

\begin{equation}  \label{Icom eq21}
(-i(1-\varphi_1)[\tilde P,{O} p^w(\lambda_1)](1-\varphi_1)u_k|u_k)\ge
C_3h_k^{-1}\|\langle x\rangle^{-s}d(x,h_kD_x)u_k\|^2-C_4N(t).
\end{equation}
From \eqref{Icom eq19}, \eqref{Incom eq20} and \eqref{Icom eq21} 
we obtain,

\begin{equation}  \label{Icom eq22}
A_1\ge C_3h_k^{-1}\|\langle x\rangle^{-s}d(x,h_kD_x)u_k\|^2-C_5N(t).
\end{equation}
Following \eqref{Icom eq16}, \eqref{Icom eq17}, \eqref{Icom eq18} and 
\eqref{Icom eq22}, we have

\begin{equation}  \label{Incom eq22bis}
N^{\prime }(t)+C_3h_k^{-1}\|\langle x\rangle^{-s}d(x,h_kD_x)u_k\|^2\le
\beta(t)+C_6N(t),
\end{equation}
where we have set $$\beta(t)= \frac{C_1}{h_k}\|\langle x\rangle ^sg_k(t)\|.\|
\langle x\rangle ^{-s}u_k(t)\|.$$
Integrating \eqref{Incom eq22bis} between $0$ and $t$ for $t\in[0,T]$ we
obtain, 
\begin{equation}  \label{Incom eq23}
N(t)+C_3h_k^{-1}\|\langle
x\rangle^{-s}d(x,h_kD_x)u_k\|^2_{L^2([0,T]\times\Omega)}\le
\int_0^T\beta(t)dt+N(0)+C_8\int_0^tN(s)ds.
\end{equation}
By Gronwall's inequality we have for $t\in [0,T]$, 
\begin{equation}  \label{Incom eq24}
N(t)\le C_7\int_0^T\beta(t)dt+ C_8N(0).
\end{equation}
Using \eqref{Incom eq24} in \eqref{Incom eq23}, we get 
\begin{align*}
\|\langle x\rangle^{-s}d(x,h_kD_x)u_k\|_{L^2([0,T]\times\Omega)}^2 &  \notag
\\
\le& {C_8}\|\langle x\rangle ^sg_k\|_{L^2([0,T]\times\Omega)}\| \langle
x\rangle
^{-s}u_k\|_{L^2([0,T]\times\Omega)}+C_9h_k\|u_k(0)\|^2_{L^2(\Omega)}.
\end{align*}
Following \eqref{eq:contradiction} and \eqref{eq:contradiction3} we obtain

\begin{equation*} 
\|\langle x\rangle^{-s}d(x,h_kD_x)u_k\|_{L^2([0,T]\times\Omega)}^2\to 0 
\text{ when }k\to+\infty.
\end{equation*}
Let $\chi (t,\tau )\in \Con_{0}^{\infty }(\R^{2})$ supported in a
neighborhood sufficiently small around $(t_{0},\tau _{0})$ and taking
account that $d$ is supported in a neighborhood of $(x_{0},\xi _{0})$, we
have

\begin{equation*}
\|\chi(t,h^2_k)d(x,h_kD_x)u_k\|_{L^2([0,T]\times\Omega)}\to 0 \text{ when }
k\to+\infty,
\end{equation*}
then $\langle \mu, \chi^2d^2\rangle=0$ thus $(x_0,t_0,\xi_O,\tau_0)\not\in
\supp \mu$.
\end{prooff}

\subsection{The microlocal defect measure vanishes on $\{a^{2}>0\}$}

The goal of this section is to prove that the microlocal defect measure
vanishes on $\{a^{2}>0\}$. More precisely we have the following proposition.

\begin{proposition}
\label{prop: vanish on a>0} Let $u_k=\psi(h^2_kP)u$ satisfying 
\begin{equation}  \label{eq1:prop a>0}
h_k^2(D_t+P)u_k-ih_ka(h_k^2P)^{1/2}(au_k)=h_kg_k,
\end{equation}
 
\begin{equation}  \label{eq2:prop a>0}
\left\Vert \left\langle x\right\rangle ^{s}g_{k}\right\Vert
_{L^{2}([0,T]\times\Omega)}^{2}+h_k\sup_{t\in \lbrack 0,T]}\left\Vert
u_{k}(t)\right\Vert _{L^{2}(\Omega)}^{2}+h_k\underset{k\rightarrow +\infty }{
\rightarrow }0,
\end{equation}
 and
\begin{equation}  \label{eq3:prop a>0}
\left\Vert \left\langle x\right\rangle ^{-s}u_{k}\right\Vert
_{L^{2}([0,T]\times\Omega)}^{2}\underset{k\rightarrow +\infty }{\rightarrow }
1.
\end{equation}
We assume that the sequence $(W_k)=(1_{[0,T]}1_\Omega u_k)$ admits a microlocal
defect measure $\mu$ then $a^2\mu=0$.
\end{proposition}

\begin{prooff}
Taking the imaginary part of the $L^{2}([0,T]\times \Omega )$ inner product
of \eqref{eq1:prop a>0} with $u_{k}/h_{k}$, we obtain, 
\begin{align}
& \Im m \lbrack
(h_{k}(D_{t}+P)u_{k}|u_{k})-i(a(h_{k}^{2}P)^{1/2}(au_{k})|u_{k})=\Im m
(g_{k}|u_{k}).  \notag \\
\end{align}
 Using that $P$ is self-adjoint, we get
\begin{align}
& \Im m (h_{k}\int_{0}^{T}\int_{\Omega }\frac{1}{2}
D_{t}|u_{k}|^{2}dxdt)-((h_{k}^{2}P)^{1/2}(au_{k})|au_{k})=\Im m(\langle
x\rangle ^{s}g_{k}|\langle x\rangle ^{-s}u_{k}).  \label{eq4=prop a>0}
\end{align}
  From \eqref{eq2:prop a>0} and \eqref{eq3:prop a>0}, we have 
\begin{equation*}
h_{k}\int_{0}^{T}\int_{\Omega }D_{t}|u_{k}|^{2}dxdt=ih_{k}\Vert
u_{k}(0)\Vert _{L^{2}(\Omega )}^{2}-ih_{k}\Vert u_{k}(T)\Vert _{L^{2}(\Omega
)}^{2}\underset{k\rightarrow +\infty }{\rightarrow }0,
\end{equation*}
and
\begin{equation*}
|(\langle x\rangle ^{s}g_{k}|\langle x\rangle ^{-s}u_{k})|\leq \Vert \langle
x\rangle ^{s}g_{k}\Vert _{L^{2}(\Omega )}\Vert \langle x\rangle
^{-s}u_{k}\Vert _{L^{2}(\Omega )}\underset{k\rightarrow +\infty }{
\rightarrow }0.
\end{equation*}
 Following \eqref{eq4=prop a>0}, we deduce 
\begin{equation}
((h_{k}^{2}P)^{1/2}(au_{k})|au_{k})\underset{k\rightarrow +\infty }{
\rightarrow }0.  \label{eq4:prop a>0}
\end{equation}
Let $\theta\in\Con^\infty_0((0,+\infty))$ with $\theta=1$ on the support of $
\psi$. Thus we have $\theta(h_k^2P)u_k=u_k$.
Let $\tilde\theta(t)=t^{-1/4}
\theta(t)$, we have $\tilde\theta\in\Con^\infty_0((0,+\infty))$ and, 
\begin{align}
(au_k|au_k)&=(a\theta^2(h^2_kP)u_k|au_k)=(a(h^2_kP)^{1/2}\tilde
\theta^2(h^2_kP)u_k|au_k)  \notag \\
&=((h^2_kP)^{1/2}\tilde\theta^2(h^2_kP)au_k|au_k)+([a,(h^2_kP)^{1/2}\tilde
\theta^2(h^2_kP)]u_k|au_k).  \label{eq5:prop a>0}
\end{align}

From Lemma 6.3 \cite{RZ}, we have 
\begin{equation}  \label{eq6:prop a>0}
\|[a,(h^2_kP)^{1/2}\tilde\theta^2(h^2_kP)]u_k\|_{L^2(\Omega)}\le
Ch_k\|u_k\|_{L^2(\Omega)}.
\end{equation}
We have also,
\begin{align}
((h^2_kP)^{1/2}\tilde\theta^2(h^2_kP)au_k|au_k)&=
\|(h^2_kP)^{1/4}\tilde\theta(h^2_kP)au_k\|^2_{L^2([0,T]\times\Omega )} 
\notag \\
&\le \|(h^2_kP)^{1/4}au_k\|^2_{L^2([0,T]\times\Omega
)}=((h^2_kP)^{1/2}au_k|au_k)\underset{k\rightarrow +\infty }{\rightarrow }0,
\label{eq7:prop a>0}
\end{align}
from \eqref{eq4:prop a>0}. Following \eqref{eq5:prop a>0}, 
\eqref{eq6:prop
a>0} and \eqref{eq7:prop a>0}, we obtain,

\begin{equation}  \label{eq8:prop a>0}
(au_k|au_k)\underset{k\rightarrow +\infty }{\rightarrow }0.
\end{equation}

Let $b(x,t,\xi ,\tau )\in \Con_{0}^{\infty }(\R^{d}\times \R\times \R
^{d}\times \R)$, we have by standard symbolic semi-classical calculus 
\begin{align}
(a^{2}(x)b(x,t,h_{k}D_{x},h_{k}^{2}D_{t})W_{k}|W_{k})=&
(b(x,t,h_{k}D_{x},h
_k^{2}D_{t})(aW_{k})|aW_{k})  \notag \\
& +h_k(r(x,t,h_{k}D_{x},h_{k}^{2}D_{t})W_{k}|W_{k}) ,  \label{eq9:prop a>0}
\end{align}
where $r(x,t,h_{k}D_{x},h_{k}^{2}D_{t})$ is bounded on 
$L^{2}([0,T]\times \R^{d})$. Thus from \eqref{eq2:prop a>0}, we have, 
\begin{equation}
h_{k}|(r(x,t,h_{k}D_{x},h_{k}^{2}D_{t})W_{k}|W_{k})|\leq Ch_{k}\Vert
W_{k}\Vert _{L^{2}([0,T]\times \R^{d})}^{2}\underset{k\rightarrow +\infty }{
\rightarrow }0.  \label{eq10:prop a>0}
\end{equation}
From \eqref{eq8:prop a>0} and using $\Vert aW_{k}\Vert _{L^{2}(\R\times \R
^{d})}^{2}=\Vert au_{k}\Vert _{L^{2}([0,T]\times \Omega )}^{2}$ we obtain, 
\begin{equation}
|(b(x,t,hD_{x},h^{2}D_{t})(aW_{k})|aW_{k})_{L^{2}(\R\times \R^{d})}|\leq
C\Vert aW_{k}\Vert _{L^{2}(\R\times \R^{d})}^{2}\underset{k\rightarrow
+\infty }{\rightarrow }0.  \label{eq11:prop a>0}
\end{equation}
 According to the definition of the microlocal defect measure $\mu $, \eqref{eq9:prop a>0}, \eqref{eq10:prop a>0} and \eqref{eq11:prop a>0} imply the Proposition~
\ref{prop: vanish on a>0}
\end{prooff}

\subsection{Propagation properties of microlocal defect measure and end of
proof}

The statement of our results requires some geometric notions which are
classical in the microlocal study of boundary problems (cf. \cite{hor} p.
424 and 430-432). \newline
Let $M=\Omega \times \mathbb{R}_{t}$. We set $$T^{*}_b M=T^{*}M\backslash\{0\}\cup T^{*}\partial M\backslash\{0\}.$$ We have the natural restriction map $$\pi:T^{*}\mathbb{R}^{d+1}_{\overline{M}}\rightarrow T^{*}_b M,$$ which is the identity on $T^{*}\mathbb{R}^{d+1}_{M}\backslash\{0\}$ (see \cite{RZ} for details).  Consider, near a point of the boundary 
$z=(x_{1},x^{\prime },t)\in \partial M$ 
a geodesic system of coordinates given by the diffeomorphism $F$ in 
(\ref{eq:U0}), for which $z=(0,0,t)$, $M=\{(x_{1},x^{\prime },t),x_{1}>0)\}$ and
the operator $D_{t}+P$ has the form (near $z$) 
\begin{equation*}
P=D_{t}+D_{x_{1}}^{2}+R(x_{1},x^{\prime },D_{x^{\prime }})+S(x,D_{x}),
\end{equation*}
with $R$ a second order tangential operator and $S$ a first order operator.
Denoting $r(x_{1},x^{\prime },\xi ^{\prime })$ the principal symbol of $R$
and $r_{0}=r|_{x_{1}=0}$, the cotangent bundle to the boundary $T^{\star
}\partial M\backslash \{0\}$ can be decomposed (in this coordinate system)
as the disjoint union of the following regions:

\begin{itemize}
\item the elliptic region $\mathcal{E}=\{(x^{\prime },t,\xi^{\prime
},\tau)\in T^{\star}\partial M\backslash\{0\};\quad r_{0}(x^{\prime
},\xi^{\prime })+\tau >0\}$,

\item the hyperbolic region $\mathcal{H}=\{ (x^{\prime },t,\xi^{\prime
},\tau)\in T^\star\partial M\backslash\{0\} ;\quad r_{0}(x^{\prime
},\xi^{\prime })+\tau <0\}$,

\item and the glancing region $\mathcal{G}=\{(x^{\prime
},t,\xi^\prime,\tau)\in T^ \star\partial M\backslash\{0\};\quad
r_{0}(x^{\prime },\xi^{\prime })+\tau =0\}$.
\end{itemize}

For the purpose of the proofs, it is important to consider the following
subsets of the glancing region:

\begin{itemize}
\item the diffractive region $\mathcal{G}_d=\{ \zeta\in \mathcal{G},
\partial_{x_{1}}r|_{x_{1}=0}(\zeta)<0\}$,

\item the gliding region $\mathcal{G}_{g}=\{ \zeta\in \mathcal{G},
\partial_{x_{n}}r|_{x_{n}=0}(\zeta)> 0\}$; we set $\mathcal{G}^2=\mathcal{G}
_d\cup\mathcal{G}_g$,

\item and $\mathcal{G}^k=\{ \zeta\in \mathcal{G},
H_{r_0}^j(\partial_{x_{1}}r|_{x_{1}=0})(\zeta)=0,\; 0\leq j<k-2,\;
H_{r_0}^{k-2}(\partial_{x_{1}}r|_{x_{1}=0})(\zeta)\neq 0\}\quad k\geq 3$,
where $$H_{r_0}=\frac{\partial r_0}{\partial \xi^{\prime }}\frac{\partial }{
\partial x^{\prime }}- \frac{\partial r_0}{\partial x^{\prime }}\frac{
\partial }{\partial \xi^{\prime }}$$.
\end{itemize}

\begin{definition}
\label{contact d'ordre infini} We say that the bicaracteristics have no
contact of infinite order with the boundary if $\displaystyle \mathcal{G}
=\bigcup _{k=2}^{+\infty }\mathcal{G}^{k}$. 
\end{definition}

Now, we recall the definition of $\nu $ the measure 
on the boundary. By the
Lemma~\ref{LemmaA6}, we see that the sequence $(1_{[0,T]}h_{k}(\frac{
\partial w_{k}}{\partial n}))$ is bounded in $L^{2}(\mathbb{R}_{t}\times
L^{2}(\partial \Omega )).$\ Therefore with the notations in (\ref{eq:plong})
\ and Proposition \ref{mesure}, we have the following Lemma.

\begin{lemma}
\label{lemma:Bord}There exists a subsequence $(W_{\sigma _{1}(k)})$ of $
(W_{\sigma (k)})$ and a Radon measure $\nu $ on $T^{\star }(\partial \Omega
\times \mathds{R}_{t})$ such that for every $b\in \Con_{0}^{\infty }(T^{\ast
}(\partial \Omega \times \mathds{R}_{t}))$ we have 
\begin{equation*}
\lim_{k\rightarrow +\infty }\left( \mathcal{O}p(b)\left( x,t,h_{\sigma
_{1}(k)}D_{x},h_{\sigma _{1}(k)}^{2}D_{t}\right) h_{\sigma _{1}(k)}\frac{1}{i
}\tfrac{\partial W_{\sigma _{1}(k)}}{\partial n},h_{\sigma _{1}(k)}\frac{1}{i
}\tfrac{\partial W_{\sigma _{1}(k)}}{\partial n}\right) _{L^{2}(\partial
\Omega \times \mathds{R}_{t})}=\left\langle \nu ,b\right\rangle .
\end{equation*}
\end{lemma}

We give now two results on propagation of support of microlocal defect
measure. The first, Proposition~\ref{Prop Propagation interieur} for point
inside $T^\star M  $ and the second, Proposition~\ref{prop:propag} at the
boundary of $M$.

\begin{proposition}
\label{Prop Propagation interieur} Let $m_{0}=(x_{0},\xi _{0},t_{0},\tau
_{0})\in T^{\star }M$ and $U_{m_{0}}$ be a neighborhood of this point in $
T^{\star }M$. Then for every $b\in \Con_{0}^{\infty }(U_{m_{0}})$, we have 
\begin{equation}
\langle \mu ,H_{p}b\rangle =0 .  \label{eq:Hp}
\end{equation}
\end{proposition}
\begin{prooff}
It is enough to prove (\ref{eq:Hp}) when $b(x,t,\xi ,\tau )=\Phi (x,\xi
)\chi (t,\tau )$ with $\pi _{x}\supp \Phi \subset V_{x_{0}}\subset \Omega $.
Let $\varphi \in \Con_{0}^{\infty }(\Omega )$ be such that $\varphi =1$ on $
V_{x_{0}}$. We introduce 
\begin{align*}
A_{k}&=\frac{i}{h_{k}}[(\Phi (x,h_{k}D_{x})\chi (t,h_{k}^{2}D_{t})\varphi
h_{k}^{2}(D_{t}+P)1_{[0,T]}w_{k},1_{[0,T]}w_{k})_{L^{2}(\Omega \times 
\mathds{R})} \\
&\quad\quad \quad -(\Phi (x,h_{k}D_{x})\chi (t,h_{k}^{2}D_{t})\varphi
1_{[0,T]}w_{k},h_{k}^{2}(D_{t}+P)1_{[0,T]}w_{k})_{L^{2}(\Omega \times 
\mathds{R})}].
\end{align*}
We claim that we have 
\begin{equation}
\lim_{k\rightarrow +\infty }A_{k}=0 .  \label{eq:damp term}
\end{equation}
We have
\begin{align}
A_{k}&=\frac{i}{h_{k}}[(\Phi (x,h_{k}D_{x})\chi (t,h_{k}^{2}D_{t})\varphi
h_{k}^{2}[D_{t},1_{[0,T]}]w_{k},1_{[0,T]}w_{k})_{L^{2}(\Omega \times 
\mathds{R})}  \notag \\
&\quad -(\Phi (x,h_{k}D_{x})\chi (t,h_{k}^{2}D_{t})\varphi
1_{[0,T]}w_{k},h_{k}^{2}[D_{t},1_{[0,T]}]w_{k})_{L^{2}(\Omega \times 
\mathds{R})}]  \notag \\
&\quad -2\Im (\Phi (x,h_{k}D_{x})\chi (t,h_{k}^{2}D_{t})\varphi
g_{k},1_{[0,T]}w_{k})_{L^{2}(\Omega \times \mathds{R})}  \notag \\
&\quad -2\Re (\Phi (x,h_{k}D_{x})\chi (t,h_{k}^{2}D_{t})\varphi
1_{[0,T]}a(h_{k}^{2}P)^{1/2}aw_{k},1_{[0,T]}w_{k})_{L^{2}(\Omega \times 
\mathds{R})}+o(1) ,  \notag
\end{align}
where we used that $(\Phi (x,h_{k}D_{x})\chi (t,h_{k}^{2}D_{t})\varphi)-(\Phi (x,h_{k}D_{x})\chi (t,h_{k}^{2}D_{t})\varphi)^*=o(1)$ by pseudo-differen\-tial calculus. 
It was proved in \cite[proof of Proposition~A.9]{RZ} that the first and the
second terms tend to zero when $k\rightarrow +\infty $. Since $
g_{k}\rightarrow 0$ in $L_{loc}^{2}$, the third term tends also to zero when 
$k\rightarrow +\infty $.\newline
For the fourth term, according to (\ref{eq9:prop a>0}) and (\ref{eq11:prop
a>0}), it is easy to see that it tends to zero. 
Thus (\ref{eq:damp term}) is proved.

In another side, it was shown in the Proposition A.9 \cite{RZ} that 
\begin{equation*}
\lim_{k\rightarrow +\infty }A_{k}=-\langle \mu ,H_{p}(\Phi \chi )\rangle .
\end{equation*}
It follows from (\ref{eq:damp term}), (\ref{eq:Hp}) that $\langle \mu
,H_{p}b\rangle =0$ if $b=\Phi \chi $, which implies our proposition.
\end{prooff}

We consider now the case of point $m_{0}=(x_{0},\xi _{0},t_{0},\tau _{0})\in
T^{\star }\mathds{R}^{d+1}$ with $x_{0}\in \partial \Omega .$ We take, as in 
\cite{RZ}, a neighborhood $U_{x_{0}}$ so small that we can perform the
diffeomorphism $F$ described in (\ref{eq:U0}).

Let $\mu $ and $\nu $ be the measures on $T^{\star }\mathds{R}^{d+1}$\ and $
T^{\star }(\partial \Omega \times \mathds{R}_{t})$\ defined in 
Proposition~\ref{mesure} and Lemma ~\ref{lemma:Bord}. 
We denote by $\tilde{\mu}$ and $
\tilde{\nu}$ the measures on 
$T^{\star }(U_{x_{0}}\times \mathds{R}_{t})$ 
and $T^{\star }(U_{x_{0}}\cap \{y_{1}=0\}\times \mathds{R}_{t})$ 
which are the pullback of $\mu $ and $\nu $ \ by the diffeomorphism 
$\tilde{F}:(x,t)\mapsto (F(x),t).$

We first recall the Lemma~A.10 established in \cite{RZ}.

\begin{lemma}
\label{lemma:A10} Let $b\in \Con_{0}^{\infty }(T^{\star }(U_{x_{0}}\times 
\mathds{R}_{t})).$ We can find $b_{j}\in \Con_{0}^{\infty }(U_{x_{0}}\times 
\mathds{R}_{t}\times \mathds{R}_{\eta ^{\prime }}^{d-1}\times \mathds{R}
_{\tau }),$ $j=0,1$\ and $b_{2}\in \Con_{0}^{\infty }(T^{\star
}(U_{x_{0}}\times \mathds{R}_{t}))$\ with compact support in $(y,t,\eta
^{\prime },\tau )$ such that with the notations of (\ref{eq:U0}),
\begin{equation*}
b(y,t,\eta ,\tau )=b_{0}(y,t,\eta ^{\prime },\tau )+b_{1}(y,t,\eta ^{\prime
},\tau )\eta _{1}+b_{2}(y,t,\eta ,\tau )(\tau +\eta _{1}^{2}+r(y,\eta
^{\prime })),
\end{equation*}
where $r$ is the principal symbol of $R(y,D^{\prime }).$
\end{lemma}

\begin{proposition}
\label{prop:propag} With the notations of Lemma~\ref{lemma:A10} for every 
$b\in \Con_{0}^{\infty }(T^{\star }(U_{0}\times \mathds{R}_{t}))$, we have 
\begin{equation*}
\langle \widetilde{\mu },H_{p}b\rangle =-\langle \widetilde{\nu },b_{1\mid
Y_{1}=0}\rangle.
\end{equation*}
\end{proposition}

\begin{prooff}
This proof is similar to the one of Proposition A.12 \cite{RZ}. We recall
some results from \cite{RZ} used to prove Proposition A.12. 
\end{prooff}

\begin{lemma}[Lemma A.13 \protect\cite{RZ}]
Let for $j=0,1$, $b_{j}=b_{j}(Y,t,\eta ^{\prime },\tau )\in \Con_{0}^{\infty
}(U_{0}\times \mathbb{R}^{d+1})$ and $\varphi \in \Con_{0}^{\infty }(U_{0})$
, $\varphi =1$ on $\pi _{Y}\supp a_{j}$. Then, 
\begin{align}
& \frac{i}{h_{k}}[((b_{0}(\Lambda _{k})+b_{1}(\Lambda
_{k})h_{k}D_{1})\varphi
h_{k}^{2}(D_{t}+P)1_{[0,T]}v_{k}|1_{[0,T]}v_{k})_{L_{+}^{2}}  \notag \\
& \quad -\int_{U_{0}^{+}}\langle (b_{0}(\Lambda _{k})+b_{1}(\Lambda
_{k})h_{k}D_{1})\varphi 1_{[0,T]}v_{k},h_{k}^{2}(D_{t}+P)\overline{
1_{[0,T]}v_{k}}\rangle dY]  \notag \\
& =-\frac{i}{h_{k}}([h_{k}^{2}(D_{t}+P),(b_{0}(\Lambda _{k})+b_{1}(\Lambda
_{k})h_{k}D_{1})\varphi 1_{[0,T]}]v_{k}|1_{[0,T]}v_{k})_{L_{+}^{2}}  \notag
\\
& \quad -(a_{1}(0,Y^{\prime },t,h_{k}D_{Y^{\prime }},h_{k}^{2}D_{t})\varphi {
_{|Y_{1}=0}}1_{[0,T]}(h_{k}D_{1}v_{k}{_{|Y_{1}=0}})|1_{[0,T]}(h_{k}D_{1}v_{k}
{_{|Y_{1}=0}}))_{L^{2}(\mathbb{R}^{d-1}\times \mathbb{R})}.  \label{Identity}
\end{align}

Here $\langle . , .\rangle $ denotes the bracket in $\mathcal{D}^{\prime }(
\mathbb{R}_{t})$.
\end{lemma}

\begin{lemma}[Lemma A.15 \protect\cite{RZ}]
\label{mu1} Let for $j=0,1,2$, $b_{j}=b_{j}(Y,t,\eta ^{\prime },\tau )\in 
\Con_{0}^{\infty }(U_{0}\times \mathbb{R}^{d})$ and $\varphi \in
 \Con_{0}^{\infty }(U_{0})$, $\varphi =1$ on $\pi _{Y}\supp b_{j}$. 
Let us set 
\begin{equation*}
L_{k}^{j}=(b_{j}(\Lambda _{k})\varphi
(h_{k}D_{1})^{j}1_{[0,T]}v_{k},1_{[0,T]}v_{k})_{L_{+}^{2}}.
\end{equation*}
Then we have for $j=0,1,2$ 
\begin{equation*}
\lim_{k\rightarrow +\infty }L_{\sigma (k)}^{j}=\langle \widetilde{\mu }
,b_{j}\eta _{1}^{j}\rangle .
\end{equation*}
\end{lemma}

The previous Lemmas still hold in our case, since they are independent of
the equation.

\begin{lemma}
\label{I^j_k,J^j_k} Let $b=b(Y,t,\eta ^{\prime },\tau )\in \Con_{0}^{\infty
}(U_{0}\times \mathbb{R}^{d+1})$ and $\varphi \in \Con_{0}^{\infty }(U_{0})$
, $\varphi =1$ on $\pi _{Y}\supp b_{j}$. For $j=0,1$ we set,
\begin{align*}
I_{k}^{j}&=(h_{k}^{-1}b(\Lambda _{k})\varphi
(h_{k}D_{1})^{j}h_{k}^{2}(D_{t}+P)1_{[0,T]}v_{k}|1_{[0,T]}v_{k})_{L_{+}^{2}},
\\
J_{k}^{j}&=\int_{U_{0}^{+}}\langle h_{k}^{-1}b(\Lambda _{k})\varphi
(h_{k}D_{1})^{j}1_{[0,T]}v_{k}|h_{k}^{2}(D_{t}+P)1_{[0,T]}v_{k}\rangle dY.
\end{align*}
Then $\lim\limits_{k\rightarrow +\infty }I_{k}^{j}=\lim\limits_{k\rightarrow
+\infty }J_{k}^{j}=0$.
\end{lemma}

\begin{prooff}
The proof is similar to the one of Lemma~A.14 \cite{RZ}. We have, 
\begin{align*}
\!\! \! \!\! I^j_k & = \frac{1}{i}[(h_k
b(\Lambda_k)\delta_{t=0}\varphi(h_kD_1)^j
v_k(0,.)|1_{[0,T]}v_k)_{L^2_+}-(h_k
b(\Lambda_k)\delta_{t=T}\varphi(h_kD_1)^j v_k(0,.)| 1_{[0,T]}v_k)_{L^2_+}] 
\notag \\
&\quad + ( b(\Lambda_k)\varphi(h_kD_1)^j 1_{[0,T]} g_k|
1_{[0,T]}v_k)_{L^2_+}+( b(\Lambda_k)\varphi(h_kD_1)^j 1_{[0,T]}
a(h_k^2P)^{1/2}av_k|1_{[0,T]}v_k)_{L^2_+}.
\end{align*}
From Lemma A.14 \cite{RZ}, the first and the second terms of the RHS in the
previous identity tend to zero.\\
Using that $\|g_k\|_{L^2}\rightarrow 0$, we can prove that the third term
tends also to zero.

Following Lemma~\ref{lemma:H} and \eqref{eq8:prop a>0} the forth term tends
to zero. 
We conclude that $I_k^j$ tends to zero. For $J_k^j$ we argue as for $I_k^j$. 
\end{prooff}

\begin{prooff}[Proof of Proposition~\protect\ref{prop:propag}]
From Proposition~\ref{prop:sopport} $(\tau +p)\mu =0$, so we have 
\begin{equation*}
\langle \widetilde{\mu },H_{p}b\rangle =\langle \widetilde{\mu }
,H_{p}(b_{0}+b_{1}\eta _{1})\rangle.
\end{equation*}
Let consider the identity (\ref{Identity}), by Lemma~\ref{I^j_k,J^j_k}, the
LHS tends to zero when $k\rightarrow +\infty $. By the semiclassical
symbolic calculus, we have 
\begin{equation*}
\frac{i}{h_{k}}[k^{2}(D_{t}+P),(b_{0}(\Lambda _{k})+b_{1}(\Lambda
_{k})h_{k}D_{1})\varphi ]=\sum_{j=0}^{2}c_{j}(\Lambda _{k})\varphi
(h_{k}D_{1})^{j},
\end{equation*}
where $c_{j}\in \Con_{0}^{\infty }(U_{0}\times \mathbb{R}^{d+1})$, $\varphi
_{1}=1$ on $\supp\varphi $, and $\{p,b_{0}+b_{1}\eta
_{1}\}=\sum\limits_{j=0}^{2}c_{j}\eta _{1}^{j}$. Hence, using Lemma~\ref{mu1}
and Lemma~\ref{lemma:Bord}, the RHS of (\ref{Identity}) tends to 
\begin{equation*}
-\langle \widetilde{\mu },H_{p}(b_{0}+b_{1}\eta _{1})\rangle -\langle 
\widetilde{\nu },b_{1}{_{|Y_{1}=0}}\rangle ,
\end{equation*}
when $k\rightarrow +\infty $.

We conclude that 
\begin{equation*}
\langle \widetilde{\mu},H_p b\rangle=\langle \widetilde{\mu}
,H_p(b_0+b_1\eta_1)\rangle=-\langle\widetilde{\nu},b_1{_{|Y_1=0}}\rangle,
\end{equation*}
which proves the Proposition~\ref{prop:propag}.
\end{prooff}

\begin{proposition}
\label{prop: mesure G} With the notations of \cite{RZ}, we have 
\begin{equation*}
\widetilde{\nu }(\mathcal{G}_{d}\cup (\bigcup\limits_{k=3}^{+\infty }
\mathcal{G}^{k}))=0.
\end{equation*}
\end{proposition}

\begin{prooff}
The proof is the same as of Lemma A.17 in \cite{RZ}.
\end{prooff}

By measure theory methods (see \cite{Bu}, \cite{B2} and \cite{RZ}), the
propagation of the measure $\mu $ along the generalized bicharacteristic
flow is equivalent to Propositions~\ref{Prop Propagation interieur}, 
\ref{prop:propag} and \ref{prop: mesure G}. 

\appendix

\section{Appendix}\label{appendix}

In this appendix, we prove some Lemmas used above.\\
We recall the Helffer-Sj\"{o}strand formula (see \cite{DA}) used extensively
in this section. To introduce it we recall some notations.\\ 
Let $\theta \in \Con_{0}^{\infty }(\R)$ and let $\varphi \in 
\Con_{0}^{\infty }(\R)$ such that $\varphi (t)=1$ if $|t|\leq 1$ 
and $\varphi
(t)=0$ if $|t|\geq 2$. Let $N\geq 2$, we set 
\begin{equation*}
\tilde{\theta}(t,\sigma )=\sum_{q=1}^{N}\frac{\theta ^{(q)}(t)}{q!}(i\sigma
)^{q}\varphi (\sigma ).
\end{equation*}
then $\tilde{\theta}\in \Con_{0}^{\infty }(\R^{2})$ and satisfies 
\begin{equation}
|\bar{\partial}\tilde{\theta}(t,\sigma )|\leq C|\sigma |^{N}\text{ where }
\bar{\partial}\tilde{\theta}(t,\sigma )=\frac{1}{2}(\partial _{t}
\tilde{\theta}+i\partial _{\sigma }\tilde{\theta})(t,\sigma ).  \label{eq:B2}
\end{equation}
We call $\tilde{\theta}$ an almost analytic extension of $\theta $. Let $P$
a self adjoint operator. We have the following Helffer-Sj\"{o}strand formula 
\begin{equation}
\displaystyle\theta (h^{2}P)=-\frac{1}{\pi }\int_{\R^{2}}\bar{\partial}
\tilde{\theta}(t,\sigma )(z-h^{2}P)^{-1}dtd\sigma \text{ where }z=t+i\sigma .
\label{eq:B3}
\end{equation}
The formula does not depend of $N$ and $\varphi$. We recall the estimates
proved in \cite{RZ}, Lemma~A.22, we have for $f=(z-h^2P)^{-1}u$ and $\Im m
z\not=0$, 
\begin{equation}  \label{eq:B4}
\|h^2Pf\|^2_{L^2(\Omega)}+\|hD_jf \|^2_{L^2(\Omega)}+\|
hV^{1/2}f\|^2_{L^2(\Omega)}+\|f \|^2_{L^2(\Omega)}\le C\frac{\langle
|z|\rangle^2}{|\Im mz|^2}\|u\|^2_{L^2(\Omega)}.
\end{equation}
Let $h_n$ a sequence such that $h_n>0$ and $h_n\to 0$ when $n\to+\infty$. In
the sequel, for simplicity we denote such a sequence by $h$. We say $h\to 0$
instead of $h_n\to 0$ when $n\to +\infty$.

\begin{lemma}
\label{lemma:A} Let $u_h$ and $g_h$ satisfying 
\begin{equation*}
\left\lbrace 
\begin{array}{l}
h^2(D_t+P)u_h-iha(h^2P)^{1/2}(au_h)=hg_h \text{ in }
[0,T]\times \Omega \\[4pt] 
u_h=0 \text{ on } [0,T]\times \partial\Omega
\end{array}
\right.
\end{equation*}
and we assume that $\|\langle x\rangle^{-s}u_h\|^2_{L^2([0,T]\times\Omega)}\le1$, 
$h\|u_h(0)\|^2_{L^2(\Omega)}\to 0$ and $\|\langle
x\rangle^sg_h\|^2_{L^2([0,T]\times\Omega)}\to 0$ when $h\to 0$. Then 
$\displaystyle \sup_{t\in[0,T]}h\| u_h(t)\|^2_{L^2(\Omega)}\to 0$.
\end{lemma}

\begin{prooff}
Let $k(t)=h\Vert u_{h}(t)\Vert _{L^{2}(\Omega )}^{2}$, using $h\partial
_{t}u_{h}=-ihPu_{h}-a(h^{2}P)^{1/2}(au_{h})+ig_{h}$, we have 
\begin{equation*}
\begin{array}{ll}
k^{\prime }(t) & =2\Re e(h\partial _{t}u_{h}(t)|u_{h}(t)) \\ 
& =2\Re e(-ihPu_{h}(t)|u_{h}(t))-2\Re e
(a(h^{2}P)^{1/2}(au_{h})(t)|u_{h}(t))+2\Re e(ig_{h}|u_{h}).
\end{array}
\end{equation*}
Using 
\begin{equation*}
\Re e(iPu_{h}(t)|u_{h}(t))=0,
\end{equation*}
and 
\begin{equation*}
\Re e(a(h^{2}P)^{1/2}(au_{h})(t)|u_{h}(t))=\Re e
((h^{2}P)^{1/2}(au_{h})(t)|au_{h}(t))\geq 0,
\end{equation*}
we obtain 
\begin{equation*}
k'(t)\le 2\|\langle x\rangle^sg_{h}(t)\Vert _{L^{2}(\Omega )}\Vert \langle x\rangle
^{-s}u_{h}(t)\Vert _{L^{2}(\Omega )}.
\end{equation*}
Thus 
\begin{equation*}
k(t)\leq k(0)+2\Vert \langle x\rangle ^{s}g_{h}\Vert _{L^{2}([0,T]\times
\Omega )}\Vert \langle x\rangle ^{-s}u_{h}\Vert _{L^{2}([0,T]\times \Omega
)}.
\end{equation*}
The assumptions and the definition of $k$ imply the Lemma.
\end{prooff}

Let $\psi : \R\to\R$ such that $\psi(t)=0$ if $t\le \alpha$ or $t\ge \beta$
where $0<\alpha<\beta$.

\begin{lemma}
\label{lemma:B} 
Let $a\in \Con_0^\infty(\R^d)$ and $s\le 1$, there exist 
$C>0 $, $h_0$ such that, if $0<h<h_0$ we have, for all 
$u\in L^2(\Omega)$, 
\begin{equation}  \label{eq:B1}
\|\langle x\rangle^s[a,\psi(h^2P)](h^2P)^{1/2}u\|^2_{L^2(\Omega)}\le
Ch^2\|u\|^2_{L^2(\Omega)}.
\end{equation}
\end{lemma}

\begin{prooff}
We prove \eqref{eq:B1} for $u\in \Con_0^\infty(\Omega)$.\\
Taking the adjoint, \eqref{eq:B1} is equivalent to 
\begin{equation*}
\Vert (h^{2}P)^{1/2}[a,\psi (h^{2}P)]\langle x\rangle ^{s}u\Vert
_{L^{2}(\Omega )}^{2}\leq Ch^{2}\Vert u\Vert _{L^{2}(\Omega )}^{2},
\end{equation*}
which is equivalent to 
\begin{equation*}
|(\langle x\rangle ^{s}[a,\psi (h^{2}P)](h^{2}\sum \partial
_{x_{j}}a_{jk}(x)\partial _{x_{k}}+h^{2}V)[a,\psi (h^{2}P)]\langle x\rangle
^{s}u|u)\leq Ch^{2}\Vert u\Vert _{L^{2}(\Omega )}^{2}.
\end{equation*}
Thus it is enough to prove 
\begin{equation}
\Vert h\partial _{x_{j}}[a,\psi (h^{2}P)]\langle x\rangle ^{s}u\Vert \leq
Ch\Vert u\Vert _{L^{2}(\Omega )} ,  \label{eq:B1.1}
\end{equation}
and 
\begin{equation}
\Vert hV^{1/2}[a,\psi (h^{2}P)]\langle x\rangle ^{s}u\Vert \leq Ch\Vert
u\Vert _{L^{2}(\Omega )}.  \label{eq:B1.2}
\end{equation}

Now we prove \eqref{eq:B1.1}. Following the Helffer-Sj\"ostrand formula,
where $\tilde\psi$ is an almost analytic extension of $\psi$, we have 
\begin{align}
h\partial_{x_j}[a,\psi(h^2P)]\langle x\rangle^s & = -\frac1\pi\int
\bar\partial\tilde\psi(z)h\partial_{x_j}[a,(z-h^2P)^{-1}]\langle
x\rangle^sdtd\sigma  \notag \\
& =\frac1\pi\int
\bar\partial\tilde\psi(z)h\partial_{x_j}(z-h^2P)^{-1}[a,z-h^2P](z-h^2P)^{-1}
\langle x\rangle^sdtd\sigma  \notag \\
& =\frac1\pi\int
\bar\partial\tilde\psi(z)h\partial_{x_j}(z-h^2P)^{-1}[a,z-h^2P]\langle
x\rangle^s(z-h^2P)^{-1}dtd\sigma+A,  \label{eq:B5}
\end{align}
where $\displaystyle A=\frac1\pi\int
\bar\partial\tilde\psi(z)h
\partial_{x_j}(z-h^2P)^{-1}[a,z-h^2P](z-h^2P)^{-1}[\langle
x\rangle^s,z-h^2P](z-h^2P)^{-1}dsd\sigma$.

We have 
\begin{equation}  \label{eq:B5.0}
[a,z-h^2P]=h^2\sum_{j=1}^d\alpha_j(x)\partial_{x_j}+h^2c(x),
\end{equation}
where $\alpha_j$ and $c$ are compact supported. Following \eqref{eq:B5}, we
have two types of terms to control.\\
First we remark that 
\begin{equation*}
(h^{2}\sum_{j=1}^{d}\alpha _{j}(x)\partial _{x_{j}}+h^{2}c(x))\langle
x\rangle ^{s}=h^{2}\beta _{j}\partial _{x_{j}}+h^{2}d(x),
\end{equation*}
where $\beta _{j}$ and $d$ are compact supported, following \eqref{eq:B5}
and estimates \eqref{eq:B4} (with $N=3$) we obtain 
\begin{equation}  \label{eq:B5.1}
\Vert h\partial _{x_{j}}(z-h^{2}P)^{-1}(h^{2}\beta _{j}\partial
_{x_{j}}+h^{2}d(x))(z-h^{2}P)^{-1}u\Vert _{L^{2}(\Omega )}\leq 
Ch\frac{\langle |z|\rangle ^{2}}{|\Im mz|^{2}}
\Vert u\Vert _{L^{2}(\Omega )}.
\end{equation}
 Thus following \eqref{eq:B2}, we have 
\begin{equation}  \label{eq:B6}
\Vert \intmod\bar{\partial}\tilde{\psi}(z)h\partial
_{x_{j}}(z-h^{2}P)^{-1}(h^{2}\beta _{j}\partial
_{x_{j}}+h^{2}d(x))(z-h^{2}P)^{-1}udtd\sigma \Vert _{L^{2}(\Omega )}\leq
Ch\Vert u\Vert _{L^{2}(\Omega )}.
\end{equation}
Second, we have 
\begin{equation*}
\lbrack \langle x\rangle ^{s},z-h^{2}P]=h^{2}\sum_{k=1}^{d}\gamma
_{k}(x)\partial _{x_{k}}+h^{2}\gamma (x),
\end{equation*}
where $|\gamma _{k}(x)|+|\gamma (x)|\leq C\langle x\rangle ^{s-1}\leq
C^{\prime }$, with the above notations, we have following \eqref{eq:B4}, 
\begin{equation}  \label{eq:B6.1}
\begin{split}
\Vert h\partial _{x_{j}}(z-h^{2}P)^{-1}(h^{2}\alpha _{j}(x)\partial
_{x_{j}}+h^{2}c(x))(z-h^{2}P)^{-1}(h^{2}\gamma _{k}(x)\partial
_{x_{k}}+h^{2}\gamma (x))(z-h^{2}P)^{-1}u\Vert \\
\leq \displaystyle Ch^{2}\frac{\langle |z|\rangle ^{3}}{|\Im mz|^{3}}\Vert
u\Vert _{L^{2}(\Omega )},
\end{split}
\end{equation}
thus, following the proof of \eqref{eq:B6}, we prove \eqref{eq:B1.1}.

To prove \eqref{eq:B1.2}, following the Helffer-Sj\"ostrand formula we have, 
\begin{equation*}
hV^{1/2}[a,\psi(h^2P)]\langle x\rangle^s=\frac1\pi\int
\bar\partial\tilde\psi(z)hV^{1/2}(z-h^2P)^{-1}[a,z-h^2P](z-h^2P)^{-1}\langle
x\rangle^sdtd\sigma.
\end{equation*}
With the notation above, it is enough to prove 
\begin{equation}  \label{eq:B7}
\|hV^{1/2}(z-h^2P)^{-1}(h^2\sum_{j=1}^d\alpha_j(x)
\partial_{x_j}+h^2c(x))(z-h^2P)^{-1}\langle x\rangle^su\|_{L^2(\Omega)}\le Ch
\frac{\langle |z|\rangle^3}{|\Im mz|^3}\|u\|_{L^2(\Omega)}.
\end{equation}
Writing $(z-h^2P)^{-1}\langle x\rangle^s=\langle
x\rangle^s(z-h^2P)^{-1}+[(z-h^2P)^{-1},\langle x\rangle^s]$, the first term
is estimated following the proof of \eqref{eq:B5.1}. To estimate the second
term, we follow the proof of \eqref{eq:B6.1}. Thus we obtain \eqref{eq:B7}
which achieve the proof of Lemma.
\end{prooff}

\begin{lemma}
\label{lemma:B bis} Let $s\in [0,1]$ and $\chi $ a smooth function such that 
$\chi=1$ for $|x|\ge 1$. We set $\chi_R(x)=\chi(x/R)$. There exists $C>0 $
such that for all $u\in L^2(\Omega)$, 
\begin{equation*}
\|(h^2P)^{1/2}\langle x\rangle^s[\psi(h^2P),\chi_R]u\|\le Ch\|u\|.
\end{equation*}
\end{lemma}

\begin{prooff}
The proof is very close to the one of Lemma~\ref{lemma:B}. By the same
argument it is sufficient to prove 
\begin{align}
\| h\partial_{x_j}\langle x\rangle^s[\psi(h^2P),\chi_R]u\|\le Ch\|u\|,
\label{Inegalite 1 lemme B bis} \\
\| hV^{1/2}\langle x\rangle^s[\psi(h^2P),\chi_R]u\|\le Ch\|u\|.
\label{Inegalite 2 lemme B bis}
\end{align}
From the Helffer-Sj\"ostrand formula, we obtain (as in \eqref{eq:B5}) 
\begin{align}
h\partial_{x_j}\langle
x\rangle^s[\psi(h^2P),\chi_R]&=\frac1\pi\int\bar\partial\tilde\psi(z)h
\partial_{x_j} (z-h^2P)^{-1}\langle
x\rangle^s[(z-h^2P),\chi_R](z-h^2P)^{-1}dtd\sigma
\label{formule 1 lemme B bis} \\
&\quad +\frac1\pi\int\bar\partial\tilde\psi(z)h\partial_{x_j}[\langle
x\rangle^s, (z-h^2P)^{-1}][(z-h^2P),\chi_R](z-h^2P)^{-1}dtd\sigma.  \notag
\end{align}
Modulo negative power of $\Im mz$, in the first term of 
\eqref{formule 1
lemme B bis} $h\partial_{x_j} (z-h^2P)^{-1}$ is bounded on $L^2(\Omega)$
and, because $\langle x\rangle^s/R$ is bounded on the support of 
$\chi^{\prime }(x/R)$, we can write 
$\langle x\rangle^s[(z-h^2P),\chi_R]$ as
a sum of term $\alpha(x)h^2\partial_{x_j}$. This yields that $\langle
x\rangle^s[(z-h^2P),\chi_R](z-h^2P)^{-1}$
is bounded on $L^2(\Omega)$ by $Ch$
modulo negative power of $\Im mz$. 
This gives the result for the first term in \eqref{formule 1 lemme B bis}.

Writing $$[\langle x\rangle^s, (z-h^2P)^{-1}]= -(z-h^2P)^{-1}[\langle
x\rangle^s, z-h^2P](z-h^2P)^{-1}$$ and arguing as for the first term, we
obtain \eqref{Inegalite 1 lemme B bis}. By the same arguments and using that 
$hV^{1/2}(z-h^2P)^{-1}$ is bounded on $L^2(\Omega)$ modulo negative power of 
$\Im mz$ (see \cite[Lemma A.22]{RZ}), we obtain 
\eqref{Inegalite 2 lemme B
bis}.
\end{prooff}

\begin{lemma}
\label{lemma:D} Let $s$ such that $|s|\leq 1$, let 
$b\in \Con^{\infty }(\overline{\Omega })$ such that 
$|b(x)|\leq C\langle x\rangle ^{s}$ and\\ $
|\partial _{x_{j}}b(x)|+
|\partial _{x_{j}x_{k}}^{2}b|\leq C\langle x\rangle^{s-1}$, 
there exist $C>0$, $h_{0}>0$ such that, if $0<h<h_{0}$ we have, for
all $u\in L^{2}(\Omega )$, 
\begin{equation*}
\Vert \langle x\rangle ^{-s}[\psi (h^{2}P),b]u\Vert _{L^{2}(\Omega )}\leq
Ch\Vert u\Vert _{L^{2}(\Omega )}.
\end{equation*}
\end{lemma}

\begin{prooff}
By Helffer-Sj\"{o}strand formula, we have, with the notation of
 Lemma~\ref{lemma:B}, 
\begin{align}
\langle x\rangle ^{-s}[\psi (h^{2}P),b] &=\frac{1}{\pi }\displaystyle\int 
\bar{\partial}\tilde{\psi}(z)\langle x\rangle
^{-s}(z-h^{2}P)^{-1}[z-h^{2}P,b](z-h^{2}P)^{-1}dtd\sigma  \label{eq:D2} \\
&=\frac{1}{\pi }\displaystyle\int \bar{\partial}\tilde{\psi}(z)\langle
x\rangle ^{-s}(z-h^{2}P)^{-1}(h^{2}\sum_{k=1}^{d}\gamma _{k}(x)\partial
_{x_{k}}+h^{2}\gamma (x))(z-h^{2}P)^{-1}dtd\sigma ,  \notag
\end{align}
where $|\gamma_k(x)|+|\gamma(x)|\le C\langle x\rangle ^{s-1}$.\\
If $s\geq 0$, following \eqref{eq:B4}, we have 
\begin{equation}  \label{eq:D3}
\Vert \langle x\rangle ^{-s}(z-h^{2}P)^{-1}(h^{2}\sum_{k=1}^{d}\gamma
_{k}(x)\partial _{x_{k}}+h^{2}\gamma (x))(z-h^{2}P)^{-1}u\Vert
_{L^{2}(\Omega )}\leq Ch\frac{\langle |z|\rangle }{|\Im mz|}\Vert u\Vert
_{L^{2}(\Omega )},
\end{equation}
thus, following the proof of \eqref{eq:B6}, we achieve the proof of Lemma in
this case.\\
If $s<0$, we write 
\begin{equation*}
\langle x\rangle^{-s}(z-h^2P)^{-1}=(z-h^2P)^{-1} \langle
x\rangle^{-s}-(z-h^2P)^{-1}[ \langle x\rangle^{-s},(z-h^2P)](z-h^2P)^{-1}.
\end{equation*}
Putting this in \eqref{eq:D2}, we obtain two terms. The first gives 
\begin{equation}  \label{eq:D4}
\|(z-h^2P)^{-1} \langle
x\rangle^{-s}(h^2\sum_{k=1}^d\gamma_k(x)
\partial_{x_k}+h^2\gamma(x))(z-h^2P)^{-1}u\|_{L^2(\Omega)}
\le C h\frac{\langle |z|\rangle}
{|\Im mz|}\|u\|_{L^2(\Omega)}.
\end{equation}
The second gives 
\begin{equation}  \label{eq:D5}
\begin{split}
\|
(z-h^2P)^{-1}\!(h^2\!\sum_{k=1}^d\tilde\gamma_k(x)
\partial_{x_k}\!\!+h^2
\tilde
f(x))(z-h^2P)^{-1}(h^2\!\sum_{k=1}^d\gamma_k(x)\partial_{x_k}\!\!+
h^2\gamma(x))(z-h^2P)^{-1}u\| \\
\le C h^2\frac{\langle |z|\rangle^2}{|\Im mz|^2}\|u\|,
\end{split}
\end{equation}
because $|\tilde\gamma_k(x)|+|\tilde\gamma(x)|\le C \langle x\rangle^{-s-1}$
. Following \eqref{eq:D4}, \eqref{eq:D5} and the Helffer-Sj\"ostrand
formula, we obtain the Lemma.
\end{prooff}

\begin{remark}
In the Lemma~\ref{lemma:D}, we can remove the assumption $|s|\le 1$, by
commuting $\langle x\rangle ^s$ with $(z-h^2P)^{-1}$ several times, but
Lemma~\ref{lemma:D} is sufficient for us in the sequel.
\end{remark}

\begin{lemma}
\label{lemma:H}Let $a\in \Con_{0}^{\infty }(\R^{d})$, there exist 
$C>0$, $h_{0}$ such that, if $0<h<h_{0}$ we have,
 for all $u\in L^{2}(\Omega )$, 
\begin{equation*}
\Vert (h^{2}P)^{1/2}a\psi (h^{2}P)u\Vert _{L^{2}(\Omega )}^{2}\leq
Ch^{2}\Vert u\Vert _{L^{2}(\Omega )}^{2}+C\Vert au\Vert _{L^{2}(\Omega
)}^{2}.
\end{equation*}
\end{lemma}

\begin{prooff}
Writing 
\begin{equation*}
(h^{2}P)^{1/2}a\psi (h^{2}P)u=(h^{2}P)^{1/2}[a,\psi
(h^{2}P)]u+(h^{2}P)^{1/2}\psi (h^{2}P)au,
\end{equation*}
then using the Lemma~\ref{lemma:B} with $s=0$,
\begin{align*}
\Vert (h^{2}P)^{1/2}a\psi (h^{2}P)u\Vert _{L^{2}(\Omega )}^{2} &\leq \Vert
(h^{2}P)^{1/2}[a,\psi (h^{2}P)]u\Vert _{L^{2}(\Omega )}^{2}+\Vert
(h^{2}P)^{1/2}\psi (h^{2}P)au\Vert _{L^{2}(\Omega )}^{2} \\
&\leq Ch^{2}\Vert u\Vert _{L^{2}(\Omega )}^{2}+C\Vert au\Vert _{L^{2}(\Omega
)}^{2},
\end{align*}
which proves the Lemma.
\end{prooff}

\begin{lemma}
\label{Borne sur x puissance s fois decomposition} For all $s\in[-1,1]$,
there exists $C>0$ such that for all $u\in \Con^\infty_0(\Omega) $and all 
$h\in(0,1]$, we have 
\begin{equation*}
\|\langle x\rangle^s\psi(h^2P)\langle x\rangle^{-s}u\|_{L^2(\Omega)}\le C \|
u\|_{L^2(\Omega)}.
\end{equation*}
\end{lemma}

\begin{prooff}
We have by Lemma~\ref{lemma:D} 
\begin{align*}
\|\langle x\rangle^s\psi(h^2P)\langle x\rangle^{-s}u\|_{L^2(\Omega)} &\le
\|\psi(h^2P)u\|_{L^2(\Omega)}+ \|\langle x\rangle^s[\psi(h^2P),\langle
x\rangle^{-s}]u\|_{L^2(\Omega)} \\
&\le \|\psi(h^2P)u\|_{L^2(\Omega)}+Ch\| u\|_{L^2(\Omega)},
\end{align*}
which  proves the Lemma.
\end{prooff}

\begin{lemma}
\label{equivalence norme H alpha} Let $\alpha\in(-1,1) $ and $s\in [-1,1]$,
then there exist $C_1>$ and $C_2>0$ such that for all 
$u\in \Con^\infty_0(\Omega) $, we have 
\begin{equation*}
C_1 \sum_{n=0}^{+\infty}h_n^{-2\alpha}\|\langle
x\rangle^s\psi(h^2_nP)u\|_{L^2(\Omega)}^2\! \le
\sum_{n=0}^{+\infty}h_n^{-2\alpha}\| \psi(h^2_nP)\langle
x\rangle^su\|_{L^2(\Omega)}^2\! \le C_2
\sum_{n=0}^{+\infty}h_n^{-2\alpha}\|\langle
x\rangle^s\psi(h^2_nP)u\|_{L^2(\Omega)}^2,
\end{equation*}
where $\psi$ was defined in Section~\ref{reduction en frequences} and 
$h_n=2^{-n}$.
\end{lemma}

\begin{prooff}
We have 
\begin{align*}
\Vert \psi (h_{n}^{2}P)\langle x\rangle ^{s}u\Vert _{L^{2}(\Omega )}^{2}&
=\Vert \psi (h_{n}^{2}P)\langle x\rangle ^{s}\sum_{k=0}^{+\infty }\psi
(h_{k}^{2}P)u\Vert _{L^{2}(\Omega )}^{2} \\
& \leq 2\Vert \psi (h_{n}^{2}P)\langle x\rangle ^{s}\sum_{k=0}^{n+1}\psi
(h_{k}^{2}P)u\Vert _{L^{2}(\Omega )}^{2} \\
& \quad +2\Vert \psi (h_{n}^{2}P)\langle x\rangle ^{s}\sum_{k=n+2}^{+\infty
}\psi (h_{k}^{2}P)u\Vert _{L^{2}(\Omega )}^{2}=2A+2B.
\end{align*}
To estimate $A$, we can write
\begin{equation*}
A\leq 2\Vert \langle x\rangle ^{s}\psi (h_{n}^{2}P)\sum_{k=0}^{n+1}\psi
(h_{k}^{2}P)u\Vert _{L^{2}(\Omega )}^{2}+2\Vert \lbrack \psi
(h_{n}^{2}P),\langle x\rangle ^{s}]\sum_{k=0}^{n+1}\psi (h_{k}^{2}P)u\Vert
_{L^{2}(\Omega )}^{2}=2A_{1}+2A_{2}.
\end{equation*}
By support properties of $\psi $ and by the Lemma~\ref{Borne sur x puissance
s fois decomposition}, we have 
\begin{equation}
A_{1}=\Vert \langle x\rangle ^{s}\psi (h_{n}^{2}P)\sum_{k=n-1}^{n+1}\psi
(h_{k}^{2}P)u\Vert _{L^{2}(\Omega )}^{2}\leq \Vert \langle x\rangle
^{s}\sum_{k=n-1}^{n+1}\psi (h_{k}^{2}P)u\Vert _{L^{2}(\Omega )}^{2}.
\label{Estimation du premier terme sur les normes alpha}
\end{equation}
By Lemma~\ref{lemma:D} we see easily that
\begin{equation*}
A_{2}\leq Ch_{n}^{2}\Vert \langle x\rangle ^{s}\sum_{k=0}^{n+1}\psi
(h_{k}^{2}P)u\Vert _{L^{2}(\Omega )}^{2}.
\end{equation*}
Summing with respect $n$, we obtain 
\begin{equation}
\sum_{n=0}^{+\infty }h_{n}^{-2\alpha }h_{n}^{2}\Vert \langle x\rangle
^{s}\sum_{k=0}^{n+1}\psi (h_{k}^{2}P)u\Vert _{L^{2}(\Omega )}^{2}\leq
\sum_{n=0}^{+\infty }\left( \sum_{k=0}^{n+1}h_{n}^{-\alpha +1}h_{k}^{\alpha
}\left( h_{k}^{-\alpha }\Vert \langle x\rangle ^{s}\psi (h_{k}^{2}P)u\Vert
_{L^{2}(\Omega )}\right) \right) ^{2}.  \label{Majoration type convolution}
\end{equation}
We have $h_{n}^{-\alpha +1}h_{k}^{\alpha }=2^{-(1-\alpha )(n-k)}2^{-k}\leq
2^{-(1-\alpha )(n-k)}$ and $(2^{-(1-\alpha )j})_{j\geq 0}\in \ell ^{1}$
because $1-\alpha >0$. We can consider the right hand side of 
\eqref{Majoration type convolution} as a convolution 
$\ell ^{1}\ast \ell
^{2}$ and we obtain the estimation of this term by 
$C\sum\limits_{n=0}^{+
\infty }h_{n}^{-2\alpha}
\Vert
 \langle x\rangle ^{s}\psi (h_{n}^{2}P)u\Vert
_{L^{2}(\Omega )}^{2}$ which estimates, with 
\eqref{Estimation du premier
terme sur les normes alpha}, the term $A$.

Now we estimate $B$. By support properties of $\psi$ and Lemma~\ref{lemma:D}
 it follows that 
\begin{align*}
B&=\| \psi(h^2_nP)\langle
x\rangle^s\sum_{k=n+2}^{+\infty}\psi(h^2_kP)
\sum_{j=k-1}^{k+1}\psi(h^2_jP)u
\|_{L^2(\Omega)}^2  \notag \\
&= \| \psi(h^2_nP)\sum_{k=n+2}^{+\infty}[\langle
x\rangle^s,\psi(h^2_kP)]\sum_{j=k-1}^{k+1}
\psi(h^2_jP)u\|_{L^2(\Omega)}^2 
\notag \\
&\le C \left( \sum_{k=n+2}^{+\infty} h_k\|\langle x\rangle^s
\sum_{j=k-1}^{k+1}\psi(h^2_jP)u\|_{L^2(\Omega)} \right) ^2 .
\end{align*}
Summing with respect $n$, we obtain
\begin{multline*}
\sum_{n=0}^{+\infty}h_n^{-2\alpha}\left( \sum_{k=n+2}^{+\infty} h_k\|\langle
x\rangle^s \sum_{j=k-1}^{k+1}\psi(h^2_jP)u\|_{L^2(\Omega)} \right) ^2 \\
\le \sum_{n=0}^{+\infty} \left( \sum_{k=n+2}^{+\infty} h_n^{-\alpha}
h_k^{1+\alpha} \left( h_k^{-\alpha} \|\langle x\rangle^s
\sum_{j=k-1}^{k+1}\psi(h^2_jP)u\|_{L^2(\Omega)} \right) \right) ^2.
\end{multline*}
We have $h_n^{-\alpha} h_k^{1+\alpha}=2^{-(1+\alpha)(k-n)}2^{-n}\le
2^{-(1+\alpha)(k-n)}$ and $(2^{-(1+\alpha)j})\in \ell^1$ since $1+\alpha>0$.
We can conclude as for the term $A$ above. We have proved the right
inequality of the Lemma.

We prove the other inequality.

We have, 
\begin{align*}
\| \langle x\rangle^s\psi(h^2_nP)u\|_{L^2(\Omega)}^2&= \| \langle
x\rangle^s\psi(h^2_nP) \langle
x\rangle^{-s}\sum_{k=0}^{+\infty}\psi(h^2_kP)\langle
x\rangle^su\|_{L^2(\Omega)}^2  \notag \\
&\le 2\| \langle x\rangle^s\psi(h^2_nP) \langle
x\rangle^{-s}\sum_{k=0}^{n+1}\psi(h^2_kP)\langle
x\rangle^su\|_{L^2(\Omega)}^2  \notag \\
&\quad +2\| \langle x\rangle^s\psi(h^2_nP) \langle
x\rangle^{-s}\sum_{k=n+2}^{+\infty}\psi(h^2_kP)\langle
x\rangle^su\|_{L^2(\Omega)}^2=2D+2E.
\end{align*}
We have by properties of support of $\psi$, 
\begin{align*}
D\le 2 \|\psi(h^2_nP)\sum_{k=n-1}^{n+1}\psi(h^2_kP)\langle
x\rangle^su\|_{L^2(\Omega)}^2+2\| [\langle x\rangle^s,\psi(h^2_nP)] \langle
x\rangle^{-s}\sum_{k=0}^{n+1}\psi(h^2_kP)\langle
x\rangle^su\|_{L^2(\Omega)}^2.
\end{align*}
The estimate of the first term is clear, for the second using 
Lemma~\ref{lemma:D}, we get
\begin{multline*}
\sum_{n=0}^{+\infty}h_n^{-2\alpha} \| [\langle x\rangle^s,\psi(h^2_nP)]
\langle x\rangle^{-s}\sum_{k=0}^{n+1}\psi(h^2_kP)\langle
x\rangle^su\|_{L^2(\Omega)}^2 \\
\le\sum_{n=0}^{+\infty} \left( \sum_{k=0}^{n+1}h_n^{-\alpha+1}h_k^\alpha
\left( h_k^{-\alpha} \| \psi(h^2_kP)\langle x\rangle^su\|_{L^2(\Omega)}
\right) \right) ^2.
\end{multline*}
We have $h_n^{-\alpha+1}h_k^\alpha\le 2^{-(1-\alpha)(n-k)}$ and we can
conclude as above by convolution argument.

For $E$, it follows from the support properties of $\psi$, Lemma~\ref{Borne sur x
puissance s fois decomposition} and Lemma~\ref{lemma:D}, 
\begin{align*}
E&=\| \langle x\rangle^s\psi(h^2_nP) \langle
x\rangle^{-s}\sum_{k=n+2}^{+\infty}\psi(h^2_kP)\sum_{j=k-1}^{k+1}
\psi(h^2_jP)\langle x\rangle^su\|_{L^2(\Omega)}^2 \\
&\le \| \langle x\rangle^s\psi(h^2_nP) \sum_{k=n+2}^{+\infty}[\langle
x\rangle^{-s}, \psi(h^2_kP)]\sum_{j=k-1}^{k+1}\psi(h^2_jP)\langle
x\rangle^su\|_{L^2(\Omega)}^2 \\
&\le C \| \langle x\rangle^s \sum_{k=n+2}^{+\infty}[\langle x\rangle^{-s},
\psi(h^2_kP)]\sum_{j=k-1}^{k+1}\psi(h^2_jP)\langle
x\rangle^su\|_{L^2(\Omega)}^2 \\
&\le C \left( \sum_{k=n+2}^{+\infty}
h_k\|\sum_{j=k-1}^{k+1}\psi(h^2_jP)\langle x\rangle^su\|_{L^2(\Omega)}
\right) ^2.
\end{align*}
Summing with respect $n$, we obtain, 
\begin{multline*}
\sum_{n=0}^{+\infty} h_n^{-2\alpha}\| \langle x\rangle^s\psi(h^2_nP) \langle
x\rangle^{-s}\sum_{k=n+2}^{+\infty}\psi(h^2_kP)\langle
x\rangle^su\|_{L^2(\Omega)}^2 \\
\le \sum_{n=0}^{+\infty} \left( \sum_{k=n+2}^{+\infty}
h_n^{-\alpha}h_k^{1+\alpha} \left(
h_k^{-\alpha}\|\sum_{j=k-1}^{k+1}\psi(h^2_jP)\langle
x\rangle^su\|_{L^2(\Omega)} \right) \right) ^2.
\end{multline*}
We have $h_n^{-\alpha}h_k^{1+\alpha}\le 2^{-(n-k)(1+\alpha)}$ and we can
conclude by convolution argument.
\end{prooff}

\begin{lemma}
\label{lemma : premier commutateur} Let $s\in[-1,1]$, $\alpha\in(-1,3/2)$
there exists $C>0$ such that for all $u\in L^2(\Omega)$, we have 
\begin{equation*}
\sum_{k=0}^{+\infty}h_k^{-1}\|\langle x\rangle^s
[\psi(h_k^2P),a](h_k^2P)^{1/2}a(h_k^2P)^{-\alpha/2}u \|_{L^2(\Omega)}^2\le
C\|u\|_{L^2(\Omega)}^2.
\end{equation*}
\end{lemma}

\begin{prooff}
Following the properties of $\psi $, we have 
\begin{equation*}
(h_{k}^{2}P)^{1/2}=\sum\limits_{j=0}^{+\infty }h_{k}h_{j}^{-1}\psi
_{0}(h_{j}^{2}P)
\end{equation*}
where $\psi _{0}(\sigma )=\sigma ^{1/2}\psi (\sigma )$ and 
\begin{equation*}
(h_{k}^{2}P)^{-\alpha /2}=\sum_{n=0}^{+\infty }h_{k}^{-\alpha }h_{n}^{\alpha
}\psi _{1}(h_{n}^{2}P)
\end{equation*}
where $\psi _{1}(\sigma )=\sigma ^{-\alpha /2}\psi (\sigma )$. Thus we must
prove, 
\begin{equation}
\sum_{k=1}^{+\infty }h_{k}^{-1}\Vert \sum_{(j,n)\in \N^{\ast
2}}h_{k}^{1-\alpha }h_{j}^{-1}h_{n}^{\alpha }\langle x\rangle ^{s}[\psi
(h_{k}^{2}P),a]\psi _{0}(h_{j}^{2}P)a\psi _{1}(h_{n}^{2}P)u\Vert
_{L^{2}(\Omega )}^{2}\leq C\Vert u\Vert _{L^{2}(\Omega )}^{2}.
\label{majoration pour le premier terme du commutateur}
\end{equation}
Let us introduce for each $k$ the following partition of $\N^{2}$. 
\begin{align*}
& A_{k}^{1}=\{(j,n)\in \N^{2},\ k\geq j-2\text{ or }k\geq n-2,
\text{ and }j\geq n-2\}, \\
& A_{k}^{2}=\{(j,n)\in \N^{2},\ k\geq j-2\text{ or }k\geq n-2,\
text{ and }j\leq n-3\}, \\
& A_{k}^{3}=\{(j,n)\in \N^{2},\ k\leq j-3\text{ and }k\leq n-3\}.
\end{align*}
In the sequel, for each set $A_{k}^{p}$ we will prove 
\eqref{majoration pour
le premier terme du commutateur}.

Let $\psi _{2}\in \Con_{0}^{\infty }(0,+\infty )$ such that $\psi _{2}=1$ on
the support of $\psi $. We have, 
\begin{equation*}
\sum_{k=0}^{+\infty }h_{k}^{-1}\Vert \sum_{(j,n)\in
A_{k}^{1}}h_{k}^{1-\alpha }h_{j}^{-1}h_{n}^{\alpha }\langle x\rangle
^{s}[\psi (h_{k}^{2}P),a]\psi _{2}(h_{j}^{2}P)\psi _{0}(h_{j}^{2}P)a\psi
_{1}(h_{n}^{2}P)u\Vert _{L^{2}(\Omega )}^{2}\leq 2A+2B,
\end{equation*}
where 
\begin{align}
A& =\sum_{k=0}^{+\infty }h_{k}^{-1}\Vert \sum_{(j,n)\in
A_{k}^{1}}h_{k}^{1-\alpha }h_{j}^{-1}h_{n}^{\alpha }\langle x\rangle
^{s}[\psi (h_{k}^{2}P),a]\psi _{2}(h_{j}^{2}P)a\psi _{0}(h_{j}^{2}P)\psi
_{1}(h_{n}^{2}P)u\Vert _{L^{2}(\Omega )}^{2}  \notag \\
& \leq C\sum_{k=0}^{+\infty }h_{k}^{-1}\left( \sum_{\substack{ (j,n)\in
A_{k}^{1}  \\ |j-n|\leq 1}}h_{k}^{1-\alpha }h_{n}^{-1+\alpha }\Vert \langle
x\rangle ^{s}[\psi (h_{k}^{2}P),a]\psi _{2}(h_{j}^{2}P)a\psi
_{0}(h_{j}^{2}P)\psi _{1}(h_{n}^{2}P)u\Vert _{L^{2}(\Omega )}\right) ^{2} 
\notag \\
& \leq C\sum_{k=0}^{+\infty }\left( \sum_{n\leq k+4}h_{k}^{3/2-\alpha
}h_{n}^{-1+\alpha }\Vert \psi _{1}(h_{n}^{2}P)u\Vert _{L^{2}(\Omega
)}\right) ^{2}\text{ (by Lemma~\ref{lemma:D})}.
\label{Inegalite premier terme A}
\end{align}
We have $h_{k}^{3/2-\alpha }h_{n}^{-1+\alpha }=2^{-(k-n)(3/2-\alpha
)}2^{-n/2}\leq 2^{-(k-n)(3/2-\alpha )}$ and we can see 
\eqref{Inegalite
premier terme A} as a convolution $\ell ^{1}\ast \ell ^{2}$ if $\alpha <3/2$
which prove \eqref{majoration pour le premier terme du commutateur} for this
term.

For $B$, we can see that
\begin{align*}
B&=\sum_{k=0}^{+\infty}h_k^{-1}\|\sum_{(j,n)\in
A^1_k}h_k^{1-\alpha}h_j^{-1}h_n^\alpha \langle x\rangle^s
[\psi(h_k^2P),a]\psi_2(h_j^2P)[\psi_0(h_j^2P),a]\psi_1(h_n^2P)u
\|_{L^2(\Omega)}^2 \\
&\le 2C+2D,
\end{align*}
where 
\begin{equation*}
C=\sum_{k=0}^{+\infty}h_k^{-1}\|\sum_{(j,n)\in
A^1_k}h_k^{1-\alpha}h_j^{-1}h_n^\alpha \langle x\rangle^s
[\psi(h_k^2P),a][\psi_0(h_j^2P),a]\psi_2(h_j^2P)\psi_1(h_n^2P)u
\|_{L^2(\Omega)}^2.
\end{equation*}
In the last sum $|j-n|\le 1$, then we can estimate this term as the term $A$.

We have 
\begin{align*}
D& =\sum_{k=0}^{+\infty }h_{k}^{-1}\Vert \sum_{(j,n)\in
A_{k}^{1}}h_{k}^{1-\alpha }h_{j}^{-1}h_{n}^{\alpha }\langle x\rangle
^{s}[\psi (h_{k}^{2}P),a][\psi _{2}(h_{j}^{2}P),[\psi
_{0}(h_{j}^{2}P),a]]\psi _{1}(h_{n}^{2}P)u\Vert _{L^{2}(\Omega )}^{2} \\
& \leq \sum_{k=0}^{+\infty }\left( \sum_{(j,n)\in
A_{k}^{1}}h_{j}h_{k}^{3/2-\alpha }h_{n}^{\alpha }\Vert \psi
_{1}(h_{n}^{2}P)u\Vert _{L^{2}(\Omega )}\right) ^{2}
\text{ (by Lemma~\ref{lemma:D} 
and Lemma~\ref{double commutateur})}.
\end{align*}
In $A_{k}^{1}$, we have $j\geq n-2$ then the sum over $j$ gives a constant
time $h_{n}$. Then, 
\begin{align*}
D& \leq C\sum_{k=0}^{+\infty }\left( \sum_{n\leq k+4}h_{k}^{3/2-\alpha
}h_{n}^{1+\alpha }\Vert \psi _{1}(h_{n}^{2}P)u\Vert _{L^{2}(\Omega )}\right)
^{2} \\
& \leq C\sum_{k=0}^{+\infty }h_{k}^{3-2\alpha }\left( \sum_{n\leq
k+4}h_{n}^{2+2\alpha }\right) \left( \sum_{n\leq k+4}\Vert \psi
_{1}(h_{n}^{2}P)u\Vert _{L^{2}(\Omega )}^{2}\right) ,
\end{align*}
by Cauchy-Schwarz inequality and as all the sums converge if $\alpha \in
(-1,3/2)$, we obtain \eqref{majoration pour le premier terme du commutateur}.

Now we will estimate the sum over $A_{k}^{2}$. We have with the 
function $\psi _{2}$ defined above, 
as $\psi_0(h_j^2P)\psi_2(h_n^2P)=0$, because $j\le n-2$,
\begin{align*}
& \sum_{k=0}^{+\infty }h_{k}^{-1}\Vert \sum_{(j,n)\in
A_{k}^{2}}h_{k}^{1-\alpha }h_{j}^{-1}h_{n}^{\alpha }\langle x\rangle
^{s}[\psi (h_{k}^{2}P),a]\psi _{0}(h_{j}^{2}P)a\psi _{2}^{2}(h_{n}^{2}P)\psi
_{1}(h_{n}^{2}P)u\Vert _{L^{2}(\Omega )}^{2} \\
& \quad =\sum_{k=0}^{+\infty }h_{k}^{-1}\Vert \sum_{(j,n)\in
A_{k}^{2}}h_{k}^{1-\alpha }h_{j}^{-1}h_{n}^{\alpha }\langle x\rangle
^{s}[\psi (h_{k}^{2}P),a]\psi _{0}(h_{j}^{2}P)[[a,\psi
_{2}(h_{n}^{2}P)],\psi _{2}(h_{n}^{2}P)]\psi _{1}(h_{n}^{2}P)u\Vert
_{L^{2}(\Omega )}^{2} \\
& \quad \leq C\sum_{k=0}^{+\infty }\left( \sum_{(j,n)\in
A_{k}^{2}}h_{k}^{3/2-\alpha }h_{j}^{-1}h_{n}^{2+\alpha }\Vert \psi
_{1}(h_{n}^{2}P)u\Vert _{L^{2}(\Omega )}\right) ^{2}
\text{ (by Lemma~\ref{lemma:D} 
and the Lemma~\ref{double commutateur}).}
\end{align*}
As $\sum\limits_{j\leq n-3}h_{j}^{-1}\leq Ch_{n}^{-1}$, we can end the
proof as for the term $D$ above.

Finally we treat the sum over $A^3_k$. We have,  
as $\psi(h^2_kP)\psi_0(h^2_jP)=0$.

\begin{align*}
& \sum_{k=0}^{+\infty }h_{k}^{-1}\Vert \sum_{(j,n)\in
A_{k}^{3}}h_{k}^{1-\alpha }h_{j}^{-1}h_{n}^{\alpha }\langle x\rangle
^{s}[\psi (h_{k}^{2}P),a]\psi _{0}(h_{j}^{2}P)a\psi _{1}(h_{n}^{2}P)u\Vert
_{L^{2}(\Omega )}^{2} \\
& \quad =\sum_{k=0}^{+\infty }h_{k}^{-1}\Vert \sum_{(j,n)\in
A_{k}^{3}}h_{k}^{1-\alpha }h_{j}^{-1}h_{n}^{\alpha }\langle x\rangle
^{s}\psi (h_{k}^{2}P)a\psi _{0}(h_{j}^{2}P)\psi _{2}(h_{j}^{2}P)a\psi
_{2}(h_{n}^{2}P)\psi _{1}(h_{n}^{2}P)u\Vert _{L^{2}(\Omega )}^{2} \\
& \leq 2E+2F,  
\end{align*}
where, 
\begin{align*}
E& \ =\sum_{k=0}^{+\infty }h_{k}^{-1}\Vert \sum_{(j,n)\in
A_{k}^{3}}h_{k}^{1-\alpha }h_{j}^{-1}h_{n}^{\alpha }\langle x\rangle
^{s}\psi (h_{k}^{2}P)[[a,\psi _{0}(h_{j}^{2}P)],\psi
_{2}(h_{j}^{2}P)][a,\psi _{2}(h_{n}^{2}P)]\psi _{1}(h_{n}^{2}P)u\Vert
_{L^{2}(\Omega )}^{2} \\
& \quad \leq \sum_{k=0}^{+\infty }\left( \sum_{(j,n)\in
A_{k}^{3}}h_{k}^{1/2-\alpha }h_{j}h_{n}^{1+\alpha }\Vert \psi
_{1}(h_{n}^{2}P)u\Vert _{L^{2}(\Omega )}\right) ^{2}.
\end{align*}
If $(j,n)\in A_{k}^{3}$, we have $j\geq k+3$ then the sum over $j$ is less
than $Ch_{k}$. We obtain, 
\begin{align*}
E& \leq C\sum_{k=0}^{+\infty }\left( \sum_{n\geq k+3}h_{k}^{3/2-\alpha
}h_{n}^{1+\alpha }\Vert \psi _{1}(h_{n}^{2}P)u\Vert _{L^{2}(\Omega )}\right)
^{2} \\
& \leq C\sum_{k=0}^{+\infty }h_{k}^{3-2\alpha }\left( \sum_{n\geq
k+3}h_{n}^{2+2\alpha }\right) \left( \sum_{n\geq k+3}\Vert \psi
_{1}(h_{n}^{2}P)u\Vert _{L^{2}(\Omega )}^{2}\right) \\
& \leq C\sum_{k=0}^{+\infty }h_{k}^{5}\Vert u\Vert _{L^{2}(\Omega )}^{2}\leq
C\Vert u\Vert _{L^{2}(\Omega )}^{2}.
\end{align*}
And we have 
\begin{align*}
F& =\sum_{k=0}^{+\infty }h_{k}^{-1}
\Vert \sum_{\substack{ {(j,n)\in A_{k}^{3}
}  \\ {|j-n|\leq 1}}}h_{k}^{1-\alpha }h_{j}^{-1}h_{n}^{\alpha }\langle
x\rangle ^{s}\psi (h_{k}^{2}P)[[a,\psi _{0}(h_{j}^{2}P)],\psi
_{2}(h_{j}^{2}P)]\psi _{2}(h_{n}^{2}P)a\psi _{1}(h_{n}^{2}P)u\Vert
_{L^{2}(\Omega )}^{2} \\
& \leq C\sum_{k=0}^{+\infty }h_{k}^{1-2\alpha }\left( \sum_{n\geq
k+3}h_{n}^{1+\alpha }\Vert \psi _{1}(h_{n}^{2}P)u\Vert _{L^{2}(\Omega
)}\right) ^{2}\text{ (by Lemma~\ref{double commutateur})} \\
& \leq C\sum_{k=0}^{+\infty }h_{k}^{1-2\alpha }\left( \sum_{n\geq
k+3}h_{n}^{2+2\alpha }\right) \left( \sum_{n\geq k+3}\Vert \psi
_{1}(h_{n}^{2}P)u\Vert _{L^{2}(\Omega )}^{2}\right) \leq
C\sum_{k=0}^{+\infty }h_{k}^{3}\Vert u\Vert _{L^{2}(\Omega )}^{2}.
\end{align*}
Which achieve the proof of Lemma.
\end{prooff}

\begin{lemma}
\label{double commutateur} Let $b\in\Con^\infty(\Omega)$ 
with support in $\{|x|\le R\}$, let $\theta_1,\ \theta_2\in 
\Con^\infty_0(\R)$, let $s\in[0,1]$ 
 there exist 
$h_0>0$ and $C>0$ such that for all $u\in L^2(\Omega)$ and 
$h\in(0,h_0)$ we have, 
\begin{equation*}
\|   \langle x\rangle ^{s}
[[\theta_1(h^2P),b],\theta_2(h^2P)]u\|_{L^2(\Omega)}\le Ch^2\|
u\|_{L^2(\Omega)}.
\end{equation*}
\end{lemma}

\begin{prooff}
We give only a sketch of proof, we use the same technic than before. By the
Helffer-Sj\"ostrand formula, we have 
\begin{equation*}
[[\theta_1 (h^{2}P),b]\theta_2(h^2P)]u=\frac{1}{\pi ^2}
\int_{\R^{4}}\bar{\partial} \tilde{\theta_1}(t_1,\sigma_1 )
\bar{\partial}\tilde{\theta_2}
(t_2,\sigma_2 )[[(z_1-h^{2}P)^{-1},b],(z_2-h^2P)^{-1}]dtd\sigma ,
\end{equation*}
where $z=(z_1,z_2)$ and $z_j=t_j+i\sigma_j$.

First, we can write
\begin{align*}
&[[(z_1-h^{2}P)^{-1},b],(z_2-h^2P)^{-1}] \\
&\quad
=(z_1-h^{2}P)^{-1}(z_2-h^2P)^{-1}[[z_1-h^{2}P,b],z_2-h^2P](z_1-h^{2}P)^{-1}(z_2-h^2P)^{-1},
\end{align*}
and $$\dsp [[z_1-h^{2}P,b],z_2-h^2P]=h^4\sum_{j,k}\gamma_{jk}(x)
\partial^2_{jk}+h^4\sum_j\gamma_j(x)\partial_j+h^4\gamma_0(x),$$ 
where the $\gamma$'s are compactly supported.
Second, as
\begin{align*}
 \langle
x\rangle ^{s}(z_1-h^{2}P)^{-1}(z_2-h^2P)^{-1}=&(z_1-h^{2}P)^{-1}(z_2-h^2P)^{-1}\langle
x\rangle ^{s}\\
&+[ \langle
x\rangle ^{s},(z_1-h^{2}P)^{-1}](z_2-h^2P)^{-1}\\
&+ (z_1-h^{2}P)^{-1}  [ \langle
x\rangle ^{s}  ,(z_2-h^2P)^{-1}],
\end{align*}
and $[ \langle
x\rangle ^{s},(z-h^{2}P)^{-1}]=-(z-h^{2}P)^{-1} [ \langle
x\rangle ^{s},   (z-h^{2}P)] (z-h^{2}P)^{-1}$, then we can obtain the Lemma by using the estimate \eqref{eq:B4} and writing the commutator $[ \langle
x\rangle ^{s}, (z-h^{2}P)]$ as in the Formula \eqref{eq:D2}.
\end{prooff}

\begin{lemma}
\label{lemme deuxieme terme} Let $s\in[-1,1]$, $\alpha<3/2$, there exists $
C>0$ such that for all $u\in L^2(\Omega)$, we have 
\begin{equation*}
\sum_{k=0}^{+\infty}h_k^{-1}\|\langle x\rangle^s
a(h_k^2P)^{1/2}[\psi(h_k^2P),a](h_k^2P)^{-\alpha/2}u \|_{L^2(\Omega)}^2\le
C\|u\|_{L^2(\Omega)}^2.
\end{equation*}
\end{lemma}

\begin{prooff}
We follow the same strategy than the one for the proof of Lemma~\ref{lemma :
premier commutateur}. We have to prove, 
\begin{equation}  \label{majoration pour le deuxieme terme du commutateur}
\sum_{k=0}^{+\infty}h_k^{-1}\|\sum_{(j,n)\in\N^2}
h_k^{1-\alpha}h_j^{-1}h_n^\alpha \langle x\rangle^s
a\psi_0(h_j^2P)[\psi(h_k^2P),a]\psi_1(h_n^2P)u \|_{L^2(\Omega)}^2\le
C\|u\|_{L^2(\Omega)}^2.
\end{equation}
If $[j-k|\ge2$ and $|n-k|\ge2$, the corresponding term in the sum is null.
If $|j-k|\le 1$ (the case $|n-k|\le1$ is symmetric and let to the reader). 
We consider two cases, the first if $n\ge k+2$, term $A$ in the sequel, and
the second if $k\ge n+2$ term $B$ in the sequel. 
\begin{align*}
A&\le C\sum_{k=0}^{+\infty} \left( 
\sum_{\genfrac{}{}{0pt}{}{|j-k|\le1}{n\ge k+2}
}h_k^{-1/2-\alpha}h_n^\alpha \| \langle x\rangle^s
a\psi_0(h_j^2P)\psi(h_k^2P)a\psi_2(h_n^2P)
\psi_1(h_n^2P)u \|_{L^2(\Omega)}
\right) ^2 \\
&\le C\sum_{k=0}^{+\infty} \left( \sum_{ n\ge
k+2}h_k^{-1/2-\alpha}h_n^\alpha \|
\psi(h_k^2P)[a,\psi_2(h_n^2P)]\psi_1(h_n^2P)u \|_{L^2(\Omega)}
 \right) ^2 \\
&\le C \sum_{k=0}^{+\infty} \left( \sum_{n\ge
k+2}h_k^{-1/2-\alpha}h_n^{1+\alpha}
 \|\psi_1(h_n^2P)u \|_{L^2(\Omega)}
\right) ^2 \\
&\le C \sum_{k=0}^{+\infty}h_k^{-1-2\alpha} \left( \sum_{n\ge k+2}
h_n^{2+2\alpha} \right) \left( \sum_{n\ge k+2}\|\psi_1(h_n^2P)u
\|_{L^2(\Omega)}^2 \right) \\
&\le C\sum_{k=0}^{+\infty}h_k \| u\|_{L^2(\Omega)}^2\le C\|
u\|_{L^2(\Omega)}^2.
\end{align*}

\begin{align*}
B&\le C\sum_{k=0}^{+\infty} \left( 
\sum_{\genfrac{}{}{0pt}{}{|j-k|\le1}{k\ge n+2}
}h_k^{-1/2-\alpha}h_n^\alpha \| \langle x\rangle^s
a\psi_0(h_j^2P)\psi_2(h_k^2P)\psi(h_k^2P)a
\psi_1(h_n^2P)u \|_{L^2(\Omega)}
\right) ^2 \\
&\le C\sum_{k=0}^{+\infty} \left( \sum_{k\ge n+2}
h_k^{-1/2-\alpha}h_n^\alpha
\| [\psi_2(h_k^2P),[\psi(h_k^2P),a]]\psi_1
(h_n^2P)u \|_{L^2(\Omega)} \right)
^2 \\
&\le C\sum_{k=0}^{+\infty} \left( 
\sum_{ k\ge n+2}h_k^{3/2-\alpha}h_n^\alpha
\| \psi_1(h_n^2P)u \|_{L^2(\Omega)} \right) ^2 \\
&\le C\sum_{k=0}^{+\infty} \left( 
\sum_{ k\ge n+2}2^{-(k-n)(3/2-\alpha)} \|
\psi_1(h_n^2P)u \|_{L^2(\Omega)} \right) ^2\le C\| u\|_{L^2(\Omega)}^2,
\end{align*}
because the last term can be seen as a convolution 
$\ell^1*\ell^2$ if $\alpha<3/2$. 
The estimations on $A$ and $B$ prove 
\eqref{majoration pour le
deuxieme terme du commutateur}.
\end{prooff}


\begin{thebibliography}{99}
\bibitem{al1} \newblock  L. Aloui. \newblock Smoothing effect for
regularized Schr\"{o}dinger equation on compact manifolds. Collect. Math., 
\textbf{59} (2008) 53--62.

\bibitem{al2} \newblock  L. Aloui. \newblock  Smoothing effect for
regularized Schr\"{o}dinger equation on bounded domains. Asymptotic
Analysis, \textbf{59} (2008) 179--193.

\bibitem{alkh} \newblock L. Aloui and M. Khenissi. \newblock Stabilisation
de l'\'{e}quation des ondes dans un domaine ext\'{e}rieur. \newblock Rev.
Math. Iberoamericana, \textbf{28} (2002) 1--16.

\bibitem{alkh2} \newblock L. Aloui and M. Khenissi. \newblock Stabilization
of Schr\"{o}dinger equation in exterior domains, \newblock Control,
Optimisation and Calculus of Variations, ESAIM, \textbf{13} (2007) 570--579.

\bibitem{alKhVo} \newblock  L. Aloui, M. Khenissi and G. Vodev. \newblock 
Smoothing effect for the regularized Schr\"odinger equation with non
controlled orbits. http://arxiv.org/abs/1201.3711v1.

\bibitem{BLR} \newblock C. Bardos, G. Lebeau and J. Rauch. \newblock Sharp
sufficient conditions for the observation, control, and stabilization of
waves from the boundary. \newblock SIAM J. Control Optimization, \textbf{30}
(1992) 1024--1065.

\bibitem{burq1} \newblock N. Burq. \newblock D\'{e}croissance de l'\'{e}
nergie locale de l'\'{e}quation des ondes pour le probl\`{e}me exterieur. 
\newblock Acta Math. \textbf{180} (1998) 1--29.

\bibitem{Bu} \newblock N. Burq. \newblock Mesures semi classiques et mesures
de d\'{e}faut, Ast\'{e}ristique, \textbf{245} (1997) 167--95.

\bibitem{B2} \newblock N. Burq. \newblock  semi-classical estimates for the
resolvant in non trapping geometries. International Mathematics research
Notices, \textbf{5} (2002) 221--41.

\bibitem{bursmot} \newblock N. Burq. \newblock Smoothing effect for Schr\"{o}
dinger boundary value problems, Duke Math. J., \textbf{123} (2004) 403--427.

\bibitem{b.g.t} \newblock N. Burq, P. G\'erard \& N. Tzvetkov. \newblock On
non linear Schr\"{o}dinger equation in exterior domain, Ann. I. H. P., 
{\large 21} (2004) 295--318.

\bibitem{co.sa1} \newblock P. Constantin \& J-C. Saut. \newblock Local
smoothing properties of dispersive equation, J. Amer. Math. Soc., \textbf{1}
(1988) 413--439.

\bibitem{doi1} \newblock S. Do\"{\i}. \newblock  Smoothing effects of 
Schr\"{o}dinger evolution groups on Riemannian manifolds, Duke Math. J., 
\textbf{82} (1996) 679--706.

\bibitem{DA} \newblock E.B. Davies. \newblock  Spectral theory and
differential operators. Cambridge studies in advanced mathematics, 
\textbf{42} Cambridge Univers. press.

\bibitem{Do} \newblock S. Do\"i. \newblock  Remarks on the Cauchy problem
for Schr\"{o}dinger type equations. Communications in Partial Differential
Equations, \textbf{21} (1996) 163--178.

\bibitem{doi05} \newblock S. Do\"i. \newblock  Smoothness of solutions for
Schr\"{o}dinger equations with unbounded potential. Publications of the
Research Institute for Mathematical Sciences, Kyoto University \textbf{41}
(2005) 175--221.

\bibitem{gerard} \textsc{P. G\'{e}rard}. \newblock Microlocal defect
measures, Com. Par. Diff. Eq., \textbf{16} (1991) 1761--1794.

\bibitem{hor} \newblock L. H\"{o}rmander. \newblock The Analysis of Linear
Partial Differential Operators. Vol III, Springer, 1985.

\bibitem{ikeh} \newblock R. Ikehata. \newblock Local energy Decay for
lineair wave equation with localized dissipation. \newblock Funkcialaj
Ekvacioj, \textbf{48} (2005) 351--366.

\bibitem{kh} \newblock M. Khenissi. \newblock \'Equation des ondes amorties
dans un domaine ext\'{e}rieur. \newblock Bull.Soc. Math. France, \textbf{131}
(2003) 211--228.

\bibitem{LaxPh} \newblock P. D. Lax and R. S. Phillips. \newblock Scattering
theory, Pure and Applied Mathematics. \newblock Academics Press, New York 
\textbf{26} 1967.

\bibitem{Leb} \newblock G. Lebeau. \newblock \'Equations des ondes amorties. 
\newblock  Algebraic and Geometric Methods in Math.Physic, A. Boutet de
Monvel and V. Marchenko (eds), Kluwer Academic, The Netherlands, (1996)
73--109.

\bibitem{Melros} \newblock R. Melrose. \newblock Singularities and energie
decay in acoustical scattering. \newblock Duke Math. J, \textbf{46} (1979)
43--59.

\bibitem{melJost} \newblock R.B. Melrose and J. Sj\"{o}strand. \newblock 
Singularities of boundary value problems I. Communications in Pure Applied
Mathematics, \textbf{31} (1978) 593--617.

\bibitem{melJost2} \newblock R.B. Melrose and J. Sj\"{o}strand. \newblock  
Singularities of boundary value problems II. Communications in Pure Applied
Mathematics, \textbf{35} (1982) 129--168.

\bibitem{MoraRalStr} \newblock C. S. Morawetz, J. Ralston and W. Strauss. 
\newblock Decay of solutions of the wave equation outside non-trapping
obstacles. \newblock Comm. Pure. Appl. Math., \textbf{30} (1977) 447--508.

\bibitem{Nak} \newblock M. Nakao. \newblock Stabilization of local energy in
an exterior domain for the wave equation with a localized dissipation. 
\newblock J. Diff. Eq. \textbf{148} (1998) 388--406.

\bibitem{rals} \newblock J. Ralston. \newblock Solution of the wave equation
with localized energy, \newblock Comm. Pure Appl. Math., \textbf{22} (1969)
807--823.

\bibitem{rauch} \newblock J. Rauch, M. Taylor. \newblock   Exponential decay
of solutions to hyperbolic equations in bounded domains. Indiana Univ Math
J, \textbf{24} (1974) 79--86.

\bibitem{RZ} \newblock  L. Robbiano, C. Zuily. \newblock  The Kato smoothing
effect for Sch\"{o}dinger equations with unbounded potentials in exterior
domains. International Mathematics Research Notices, (2009) 1636--1698.

\bibitem{sjolin} \newblock  P. Sj\"{o}lin. \newblock Regularity of solutions
to the Schr\"{o}dinger equation. Duke Math. J. \textbf{\ 55} (1987) 699--715.

\bibitem{vainb} \newblock  B. Vainberg. \newblock Asymptotic methods in
equations of mathematical physics, Gordon \& Breach Science Publishers, New
York, 1989.

\bibitem{vega} \newblock  L. Vega. \newblock Schr\"{o}dinger equations:
pointwise convergence to the initial data. Proc. Amer. Math. Soc. \textbf{102
} (1988) 874--878.

\bibitem{V.Z} \newblock A. Vazy and M. Zworski. \newblock  Semiclassical
Estimates in Asymptotically Euclidean Scattering. Commun. Math. Phys., 
\textbf{212} (2000) 205--217.

\bibitem{wilcox} \newblock C. H. Wilcox. \newblock  Scattering Theory for
the d'Alembert Equation in Exterior Domains. Lecture Notes, Math \textbf{442}
Springer, New York, 1975.
\end{thebibliography}
\end{document}